\magnification=\magstep1
\input amssym.def
\input amssym
\font\ams=msbm10

\def\ctr{\centerline}
\def\ga{{\goth a}}
\def\gb{{\goth b}}
\def\gg{{\goth g}}
\def\gh{{\goth h}}
\def\gl{{\goth l}}
\def\gm{{\goth m}}
\def\gn{{\goth n}}
\def\gp{{\goth p}}
\def\gs{{\goth s}}
\def\bBo{{\bf B}}
\def\bG{{\bf G}}
\def\bH{{\bf H}}
\def\bL{{\bf L}}
\def\bM{{\bf M}}

\def\bP{{\bf P}}
\def\bS{{\bf S}}
\def\bT{{\bf T}}
\def\bV{{\bf V}}

\def\bO{{\bf O}}
\def\a{\alpha}
\def\b{\beta}
\def\g{\gamma}
\def\d{\delta}
\def\l{\lambda}
\def\s{\varsigma}

\def\Nscr{{\cal N}}
\def\Fscr{{\cal F}}

\def\Oscr{{\cal O}}

\def\Vscr{{\cal V}}

\def\Iscr{{\cal I}}

\def\Wscr{{\cal W}}
\def\Xscr{{\cal X}}

\def\Sscr{{\cal S}}
\def\M{{\cal M}}
\def\Mscr{{\cal M}}
\def\Lie{{\rm Lie}\,}
\def\sh{{\rm sh}\,}

\def\rk{{\rm rk}\,}
\def\gr{{\rm gr}\,}
\def\Det{{\rm Det}\,}

\def\Id{{\rm Id}\,}

\def\Ann {{\rm Ann}\,}

\def\End{{\rm End}\,}
\def\wt{{\rm wt}\,}
\def\ch{{\rm ch}\,}
\def\Max{{\rm Max}\,}
\def\sign{{\rm sign}\,}

\def\qed{\hbox{\hskip 4pt
                \vrule width 5pt height 6pt depth 1.5pt\hskip 2pt}}
\def\QED{\par\hfill\qed \par}
\def\pr{^{\prime}}
\def\prpr{^{\prime\prime}}
\def\st{\subset}
\def\sp{\supset}
\def\ra{\rightarrow}
\def\pmat{\pmatrix}
\def\Co{\Bbb C} 
\def\Na{{\Bbb N}}  
\font\mas=msbm10
\def\zline{\hbox{\mas\char'132}} 
\def\Pf{\parno{\bf Proof.}\par}

\def\sr{\scriptscriptstyle}

\def\vb{\vrule height 14pt depth 7pt} 
\def\ts{\tabskip0pt}	
\def\vs{\noalign{\vskip-8pt}}
\def\ss{\noalign{\vskip -1.5pt}}
\def\spro{\hbox{\ams\char'157}}
\def\ov{\overline}
\def\ua{\uparrow}
\def\da{\downarrow}
\def\la{\leftarrow}
\def\raeq{\buildrel \sim \over \ra}
\def\parno{\par\noindent}
\def\disp{\displaystyle}
\voffset=-2.2truecm
\hsize=16truecm
\vsize=24truecm
\topskip=2.5 truecm
\hfuzz 5pt
\baselineskip=18truept

\tolerance 1000
\topglue 5truecm

\centerline{{\bf Quantization of hypersurface orbital varieties 
in $\gs\gl_n$}\footnote {}{Supported in part by the Minerva Foundation, Germany,
Grant No. 8466}}
\vskip 1truecm
\centerline{\bf Anthony Joseph}
\ctr{The Donald Frey Professorial Chair}
           \ctr{Department of Mathematics}
           \ctr{The Weizmann Institute of Science}
           \ctr{Rehovot 76100, \ \ Israel}
\ctr{and}
\ctr{ Institut de Math{\'e}matiques (UMR 7586)}
        \ctr{ 175 rue du Chevaleret, Plateau 7D}
        \ctr{Paris 75013 Cedex, \ \ France}
           \ctr{{\it e-mail:} joseph@wisdom.weizmann.ac.il} 

\bigskip
\centerline{and}
\bigskip
\centerline {\bf Anna Melnikov}
\centerline {Department of Mathematics}
\centerline {Haifa University}
\centerline {Mount Carmel}
\centerline {Haifa 31905, Israel}
\ctr{{\it e-mail:} anna@wisdom.weizmann.ac.il}
\vskip 2truecm 
\parno
{\bf 1. \ \ Introduction.}
\parno
{\bf 1.1.} \ \ Amongst his many great contributions to mathematics, 
Alexander Kirillov was a co-founder of the ``orbit method''. This exploits the 
symplectic structure on coadjoint orbits with respect to a Lie algebra $\gg.$
From a Lagrangian subvariety one attempts to construct a representation of
$\gg$ associated to the given orbit. If $\gg$ is semisimple, the orbit
is nilpotent and one wishes to construct a highest weight module, then this 
Lagrangian subvariety should be a so-called orbital variety. Then ``quantization'' 
of the latter leads to the required representation [J3].
\par
In this paper we consider only the case when $\gg=\gs\gl_n(\Co)$. By [M4] all
orbital varieties can be weakly quantized (though this fails [J3, 1.3] for
arbitrary $\gg$). It remains to show that they can be strongly quantized.
This is a rather more delicate question. We shall settle this positively
for orbital varieties which are of codimension 1 in the nilradical of a parabolic.
These are called hypersurface orbital varieties.
\parno
{\bf 1.2.} \ \ Let $\gg=\gn^-\oplus\gh\oplus\gn^+$ be a triangular decomposition.
We identify $\gn^-$ with $(\gn^+)^*$ through the Killing form. Let $\bG$ be
adjoint group of $\gg$ and $\bBo$ the Borel subgroup with Lie algebra\ $\gb:=\gh\oplus
\gn^+$. For any Lie algebra\ $\ga$, let $S(\ga)$ (resp. $U(\ga)$) denote its symmetric
(resp. enveloping) algebra. It is well-known that $\bG\gn^+$ is a finite union (of so-called
nilpotent) orbits classified by Jordan normal form. Let $\Oscr$ be such an orbit.
After Spaltenstein [S] the irreducible\ components of $\Oscr\cap\gn^+$ have dimension
$1/2 \dim\Oscr$. They are called orbital varieties and after Steinberg [St]
are classified by  the standard tableaux whose shape is specified by
$\Oscr$. As noted in [J1, 7.3] an orbital variety $\Vscr$ can be characterized
as an irreducible\ subvariety of $\gn^+$ for which the ideal $I(\ov\Vscr)\st S(\gn^-)$ 
of definition of its closure is stable under the Poisson bracket on
$\gg$ (induced by the Lie bracket).
\parno
{\bf 1.3.} \ \ Retain the notation of 1.2. Roughly speaking $\Vscr$ can be
strongly quantized if $S(\gn^-)/I(\ov\Vscr)$ can be given the structure of a highest weight module\
(not necessarily simple). Here the choice of highest weight\ is a crucial and delicate
point. Take $\l\in\gh^*$ and let $V(\l)$ be a highest weight module\ with highest weight\ vector $v_\l$
of weight $\l$. Let $\Fscr$ be the canonical (degree) filtration on $U(\gn^-)$ 
(resp. on $S(\gn^-)$). We say that $V(\l)$ is a strong (resp. weak)
quantization\ of $\Vscr$ if $\gr_{\Fscr}\Ann_{U(\gn^-)}v_\l=I(\ov\Vscr)$ (resp. if its radical
equals $I(\ov\Vscr)$\ ).
\parno
{\bf 1.4.} \ \ In order to exhibit a strong quantization\ we have to compute $I(\ov\Vscr)$ (or
at least to calculate the formal character of $S(\gn^-)/I(\ov\Vscr)$\ ). In [J2, Lecture 7]
the general form that this should take is suggested. Let $\Nscr$ be the
generic matrix of $\gn^-+1d$. Then $I(\ov\Vscr)$ should take the form $\gr_{\Fscr}<a_1,
a_2,\cdots,a_n>$, where $a_i$ are amongst the minors of $\Nscr$. This suggestion derived from an algorithm for 
$\Ann_{U(\gn^-)}\Vscr$ based on the Enright functor [J1, 8.4] together with calculations in the thesis
of E. Benlolo.
\parno
{\bf 1.5.} \ \ Of course the above does not tackle the difficult question as
to which minors $a_i$ to choose. Recently Benlolo and Sanderson [BS] 
made an appealing conjecture for 
this choice concerning orbital varieties $\Vscr$ having codimension 1
in the nilradical $\gm^+$ of a parabolic subalgebra\ $\gp$. Apart from the obvious
$1\times 1$ minors it follows from  Krull's theorem 
that $I(\ov\Vscr)$ is generated by just
one element $f$ which can be assumed to be an $\gh$ weight vector. A result in
[M3] describes in particular\ all the orbital varieties of codimension 1 in a given $\gm^+$.
Using this, Benlolo-Sanderson [BS] show that one may reduce to the case when
$f$ cannot be expressed in terms of the generic matrices of $\gn^-$ obtained
by deleting the first columns and last rows. Up to this reduction they show that
$\gm^+$ is given by blocks of size $c_{\sr 1},c_{\sr 2},\cdots
,c_l$ with $c_{\sr 1}=c_l=c$ and $c_i\neq c,\ \forall\ i:1< i< l$. Now let
$\Mscr(t)$ be the generic matrix of $\gm^-+t1d$, where $\gm^-$ is the opposed
algebra of $\gm^+$ identified with $(\gm^+)^*$  
and $t$ is an indeterminate. Let $M(t)$ be the $(n-c)\times(n-c)$
minor of $\Mscr(t)$ lying in the bottom left-hand corner. Then $M(t)$ is a
polynomial\ in $t$. They conjecture that $f$ is the coefficient\ of its lowest degree term.
Let $\b_1$ be the highest root of $\gn^+,\ \b_2$ the highest root of the subalgebra\
of $\gn^+$ obtained by omitting the first row and last column, and so on.
They show that $f$ has weight $-\disp{\sum^c_{i=1}}\b_i$.
\parno
{\bf 1.6.} \ \ One may note that  $f$ of 1.5 can also be described as
$\gr_{\Fscr}M(1)$. However as we shall see there is a certain computational
advantage in retaining $t$.  (Furthermore BS compute the lowest power of $t$
in $M(t)$). One may also replace $\gm^-+1d$ in the above by $\gn^-+1d$ as
long as one includes the obvious $1\times 1$ minors in $I(\ov\Vscr)$. In this sense
$M(1)$ only depends on the first column size $c$. Likewise $f=\gr_{\Fscr}M(1)$
depends only partially on $\gm^-$. In particular the weight of $f$ only depends 
on $c$. Although
BS do check their conjecture in a number of cases, it is interesting to
know if it holds in general since the above independence\ is somewhat surprising.
Moreover even the relatively simple case of hypersurface varieties has some 
remarkable structure.
\parno
{\bf 1.7.} \ \ Our proof of the BS conjecture involves showing that 
$S(\gm^-)f$ is semi-prime and its zero variety is  $\bP$ stable, where $\bP\sp \bBo$
is the parabolic subgroup with Lie algebra $\gp$. This holds irrespective of whether
$c_i\neq c$, for $i\neq 1,l$. However we show that $f$ is irreducible\ if and only if\ this
condition holds. This is a delicate point. We shall give three proofs. One is
based on a careful counting of $t$ powers in multiplying out $M(t)$ and a
knowledge [M2, M3] 
of the orbital varieties of codimension 1 in $\gm^+$. The second is based on a
representation theoretic argument combined with [M4] which asserts that the
associated variety of a highest weight module\ of {\it integral} highest weight\ is irreducible. 
The third is the least computational and is based on the irreducibility
of the associated variety of the annihilator of a simple highest weight 
module.
\parno
{\bf 1.8.} \ \ We shall exhibit a strong quantization\ of $\Vscr$ by determining a
simple highest weight module\ $V(\l)$ whose formal character $\ch V(\l)$ coincides with that of
$S(\gm^-)/S(\gm^-)f$ and by using the linear independence\ of the characteristic polynomials
of orbital varieties. To determine $\l$ and compute $\ch V(\l)$, we use the formula
of Jantzen [Ja] describing the Shapovalov determinants defined with respect to
$\gp$. Let $\Iscr$ be the subset of the set of simple roots $\Pi$ defining
$\gp$. Assume $(\a,\a)=2$ for any root. Then it suffices to take $\l\in\gh^*$
such that $(\l,\a)=0, \ \forall\ \a\in\Iscr, \ (\l,\b_1)=c-(n-1)$
 and $(\l,\g)\not\in
\zline$, for any proper sum of roots in $\Pi\pr:=\Pi\setminus\Iscr$. However
we also give a choice of integral $\l$. This is more delicate since the
condition $c_i\neq c,\ \forall\ i:1\leq i\leq k-1$ must be invoked. Combined
with [M4] it leads to the second proof of the irreducibility of $f$.
The non-integral case gives our third proof. 
\par
We further show that there is a choice of $\l$ such that the
annihilator of $V(\l)$ is maximal as well as  $V(\l)$ being a strong
quantization of $\Vscr.$

\parno
{\bf 1.9.} \ \ The truth of BS conjecture has the following remarkable
consequence.  One may specify positive integers $s,t$ such that $f$
is exactly the highest common divisor of the $t\times t$ minors 
of $\Mscr
(0)^s$. There seems to be no elementary proof of this purely combinatorial
fact.
\bigskip
\parno
{\bf 2. Combinatorial Preliminaries and the Benlolo-Sanderson Conjecture}
\parno
{\bf 2.1}\ \ \ 
The base field is assumed algebraically closed of characteristic
zero and can be taken to be the complex field $\Co$ without loss
of generality.
\par
Recall the notation and the hypotheses of 1.1 and 1.2. 
Put $\gn=\gn^+.$
In particular an orbital variety $\Vscr$ associated to a nilpotent
orbit $\Oscr$ is an irreducible component of $\Oscr\cap \gn.$
Let $\bV$ denote the set of all orbital varieties of $\gg$ and 
$W$ the Weyl group for the pair $(\gg,\ \gh).$ We describe
first the Steinberg map of $W$ onto $\bV.$
Let $R\st\gh^*$ denote the set of non-zero roots,  $R^+$ the set
of positive roots corresponding to $\gn$ in the triangular 
decomposition of $\gg$ and $\Pi\st R^+$ the resulting set of 
simple roots. Let $X_\a={\Co} x_{\a}$ denote the root subspace 
corresponding to $\a\in R.$ Then $\gn=\bigoplus\limits_{\a\in R^+}X_{\a}$
(resp.$\gn^-=\bigoplus\limits_{\a\in -R^+}X_{\a}$).
\par
For each $w\in W$ set 
$\gn\cap^w\gn:=\bigoplus\limits_{\a\in R^+\cap w(R^+)}X_{\a}.$
For each closed, irreducible  subgroup  $\bH$ of $\bG$ let
$\bH(\gn\cap^w\gn)$ be the set of $\bH$
conjugates of $\gn\cap^w\gn.$ It is an irreducible locally 
closed subvariety. Since there are only finitely many nilpotent 
orbits in $\gg$ it follows that there exists a unique nilpotent orbit
$\Oscr$ such that $\ov{\bG(\gn\cap^w\gn)}=\ov\Oscr.$ 
A result of Steinberg [St] asserts that 
$\Vscr_w:=\ov{\bBo(\gn\cap^w\gn)}\cap\Oscr$ is an orbital
variety and that the map
$\phi:w\mapsto \Vscr_w$
is a surjection of $W$ onto $\bV.$
\parno
{\bf 2.2}\ \ \  
It is convenient to replace $\gs\gl_n$ 
by $\gg=\gg\gl_n.$ This obviously makes no difference.
Note that the adjoint action is just conjugation by $\bG=GL_n.$
Let $\gn$ be  the subalgebra
of strictly upper-triangular matrices and let $\gn^-$ 
be the subalgebra of strictly lower-triangular matrices.
Let $\bBo$ be the (Borel) subgroup of upper-triangular
matrices in $\bG.$ All parabolic
subgroups we consider further are standard, that is contain
$\bBo.$
\par
Let $e_{i,j}$ be
the matrix having $1$ in the $ij-$th entry
and $0$ elsewhere. 
\par
Take $i<j$ and let $\a_{i,j}$ be the root corresponding to
$e_{i,j}.$ Set $\a_{j,i}=-\a_{i,j}.$ We write $\a_{i,i+1}$
simply as $\a_i.$ Then $\Pi=\{\a_i\}_{i=1}^{n-1}.$ Moreover
$\a_{i,j}\in R^+$ exactly when $i<j.$ For each $\a\in \Pi,$
let $s_{\a}\in W$ be the corresponding reflection and set 
$s_i=s_{\a_i}.$ 
\parno
{\bf 2.3}\ \ \ 
We represent every element of the symmetric
group $\bS_n$  in word form 
$$ w =[a_{\sr 1}, a_{\sr 2},\ldots, a_n]\ ,\quad {\rm where}
                                    \ a_i=w(i).$$
We identify $W$ with $\bS_n$ by taking $s_i$ to be the 
elementary permutation interchanging $i,\ i+1.$
\proclaim Definition. Given $w=[a_1,\cdots,a_n].$ 
Set $p_w(i)=j$ if $a_j=i,$ that is
$p_w(i)$ is the place (index) of $i$ in the word form
of $w.$
\par
One has $w(p_w(i))=w(j)=a_j=i,$ that is $p_w(i)=w^{-1}(i).$
On the other hand $w(\a_{i,j})=\a_{w(i),w(j)}.$ Set 
$S(w)=\{\a\in R^+\ |\ w(\a)\in -R^+\}.$
\par 
This gives the following result.
\proclaim Lemma. Take $i<j.$ Then $\a_{i,j}\in S(w^{-1})$
if and only if $p_w(i)>p_w(j).$
\par
\parno
{\bf Remark.}\ \  In particular $\a_i\in S(w^{-1})$ exactly
when $i+1$ comes before $i$ in the word form of $w.$
This is of course well-known.
\parno
{\bf 2.4}\ \ \
Let $P(n)$ denote the set of partitions of $n.$ Then 
$\l:=\{\l_1\geq\l_2\geq\cdots\geq \l_k>0\}\in P(n)$
defines the sizes $\l_i$ of the Jordan blocks of some element
$x_\l\in\gn.$ Set $\Oscr_\l=\bG x_\l.$ Let 
$\l^*=\{\l^*_1\geq \l^*_2\geq \cdots\geq \l^*_l>0\}$
be the dual partition. By definition 
$\l^*_i=\sharp\{j\ |\ \l_j\geq i\}.$ Then $k=\l^*_1,\ l=\l_1$
and this notation will be fixed throughout. We view $\l$
as a Young diagram $D_\l$ with $k$ rows of length 
$\l_1,\ldots,\l_k.$ 
\par
The rank of a Jordan block of size $r$ is just $r-1$ so
$\rk x_{\l}=\sum_{i=1}^k (\l_i-1)=n-\l^*_1.$
\par
For all integer $i\geq 1,$ the Jordan blocks of $x_\l^i$
are obtained by deleting the first $i$ columns of $\l.$
Then
\item {(1)}\ \ \ $\rk x^i_\l=n-\sum_{j=1}^i\l^*_j.$
\item {(2)} Again (cf. [H, \S 3.8])
$\dim \Oscr_{\l}=n^2-\sum_{i=1}^k (\l_i^*)^2.$
\parno
{\bf 2.5}\ \ \ 
Define a partial order on $P(n)$ as
follows. Let $\l=(\l_1\geq\l_2\geq\cdots \l_k)$ and
$\mu=(\mu_1\geq\mu_2\geq\cdots \mu_k)$ be partitions
of $n$ which we can assume to correspond to diagrams of the same height by adding empty 
rows. Set $\l\leq \mu$ if
$$ \sum^i_{j=1}\l_j \geq \sum^i_{j=1}\mu_j \ ,
\qquad {\rm for \quad all}\qquad i=1,2,\cdots,k.$$
The following result of M. Gerstenhaber (cf. [H, \S 3.10]) shows
that this order corresponds to inclusion of nilpotent orbit closures:
\proclaim Theorem. Given two partitions $\l$ and $\mu$ of $n$
one has $\l\leq\mu$ if and only if $\ov\Oscr_{\mu}\st\ov\Oscr_{\l}.$
\par
In our convention $\l=(n)$ is the minimal partition of $n$ and
$\l=(1,1,\ldots,1)$ is the maximal one.
\parno
{\bf 2.6}\ \ \ 
For a partition $\l$ of $n$ one can
fill the 
boxes of $D_{\l}$ with $n$ distinct positive integers.
If the entries increase in rows from left to right
and in columns from top to bottom, we call such an array
a Young tableau. If the numbers in Young tableau
form a set of integers from $1$ to $n$ we call it standard. 
Let $\bT_n$ be the set of standard Young tableaux
of size $n.$
\par
We shall not distinguish between Young and standard tableaux. 
Indeed if a Young tableau has entries $t_1<t_2<\cdots<t_n$ then
it can be identified with the standard tableau obtained by
replacing $t_i$ by $i.$ Similarly the sequence 
$[t_{w(1)},\ldots,t_{w(n)}]:\ w\in \bS_n$ can be viewed as
a word form of $w$ by replacing $t_i$ by $i.$
We call this process standardization. Occasionally it can
cause confusion; but in such cases adequate warning will be given.
We may also concatenate word forms so then $[a_1,a_2,\ldots,a_n]$
can be written as $[a\pr,a\prpr]$ where $a\pr=[a_1,a_2,\ldots,a_i],\ 
a\prpr=[a_{i+1},\ldots,a_n].$
\par
The shape of a Young tableau $T$ is defined
to be the Young diagram from which $T$ was built.
It defines a partition of $n$
which we denote by $\sh T.$  
\par 
Given $T\in \bT_n,$ let $\sigma(T)$ denote the conjugate standard 
tableau obtained by rotation about the main diagonal. We remark
that $\sh\sigma(T)=(\sh T)^*.$ Again $\sigma$ takes a row tableau
$R$ into a column tableau $C$ and vise-versa.
\par
Given any sequence $S$ of strictly increasing integers 
(for example a row $R$ or column $C$ above) let $\widehat S$
denote the reversed sequence of strictly decreasing integers.
For any tableau or word form $S,$ let $<S>$ denote the set of entries of $S$
and $|S|$ the cardinality of $<S>.$
\par
The Robinson - Schensted correspondence
$w\mapsto (Q(w),R(w))$
gives a bijection (see, for example [Kn] or [F]) from the 
symmetric group
$\bS_n$ onto the  pairs of standard
Young tableaux of the  same shape.
By R. Steinberg [St] for all
$w,y\in \bS_n$ one has $\Vscr_w=\Vscr_y$ iff
$Q(w)=Q(y).$ This parameterizes the set of orbital
varieties $\bV$ by $\bT_n.$
Moreover $\sh Q(w)=\l$ 
if and only if $\Vscr_w$ is contained in 
$\Oscr_{\lambda}.$
\par
If $Q(w)=T$ we set $\Vscr_T:=\Vscr_w$ and $T_{\Vscr_w}:=Q(w).$ 
\par
We define a partial order on Young tableaux and a partial
order on orbital varieties analogous to the order on partitions
given in 2.5. It is called the geometric order on $\bV$ or on $\bT_n.$ 
\proclaim Definition. (i)\ \ Given distinct 
$\Vscr_1, \Vscr_2\in \bV.$ Set $\Vscr_2>\Vscr_1$ if $\Vscr_2
\subset\ov\Vscr_1.$
\item{(ii)} Given distinct $T_1, T_2\in \bT_n.$
Set $T_2>T_1$ if $\Vscr_{T_2}>\Vscr_{T_1}.$
\par
In general the problem of the combinatorial description of geometric order is extremely
difficult. Partial results are described in [M5], which  
for the case studied below
gives a complete answer. All that is needed here in 2.17; 
though a more general result is given in [M2].
\parno  
{\bf 2.7}\ \ \ 
Fix $T\in \bT_n$ with $\sh T=\l.$
For all $i,j\in \{1,\ldots,n\}$ let $T_j^i$ denote the 
$ij$-th entry of $T$ when it is defined.
\proclaim Definition. For any tableau $T$
\item{(i)} Given $s=T_j^i,$ set $r_{\sr T}(s)=i$ and 
$c_{\sr T}(s)=j.$
\item{(ii)} For all $1\leq i\leq\l^*_1=k,$ set 
$T^i=(T^i_1,\cdots, T^i_{\l_i}).$ It is the $i-$th row of $T.$
\item {(iii)} For all $j\ :1\leq j\leq \l_1=l,$ set 
$T_j=(T^1_j,\ldots,T_j^{\l_j^*}).$ It is the $j-$th column of $T.$ 
\item {(iv)} The hook number of the $ij-$th entry of $T$ is 
defined by 
$$h(T^i_j)= 1+(\l_j^*-j)+(\l_i-i).$$
\item {(v)} $T_j^i$ is called a corner entry of $T$ if $h(T^i_j)=1.$
\item{(vi)} If $i\leq j,$ then $T^{i,j}$ denotes the subtableau 
with rows $T^i,T^{i+1},\ldots,T^j.$ Otherwise it denotes the 
empty tableau.
\par
A row (resp. column) of $T$ is determined by its entries since
these must increase from left to right (resp. from top to bottom).
\par
Let $T=(T_1,T_2,\ldots,T_l),\ S=(S_1,S_2,\ldots,S_{l\pr})$
be Young tableaux given by their columns. Assume that $T,S$ have no
common entries. Then we define $(T,S)$ to be the array
whose rows  are the same as the rows
of $(T_1,T_2,\ldots,T_l,S_1,S_2,\ldots,S_{l\pr}),$ that is
$r_{\sr(T,S)}(T^i_j)=i$ and $r_{\sr (T,S)}(S_j^i)=i,$ and
ordered in the increasing order.
Of course this involves the shuffling of numbered boxes within a 
row.  
\proclaim Lemma. $(T,S)$ is a Young tableau. 
\par
\Pf
The proof is by induction on rows.
Suppose that the integers in the first $r$ rows have been placed in
increasing order and that the resulting array
of $r$ rows obtained from $(T,S)$ is a Young tableau.
Let $i_1<i_2<\cdots<i_n$ be the integers of the $r-$th row and 
$j_1<j_2<\cdots<j_m\ :\ m\leq n$
the ordered set of entries of $<T^{r+1},S^{r+1}>.$
Since $T,S$ are Young tableaux and their entries are distinct from
one another, there exists an injective map 
$\phi:\{1,2,\ldots,m\}\ra \{1,2,\ldots,n\}$ such that $j_t>i_{\phi(t)}.$
Then $j_s>j_t>i_{\phi(t)}$ for all $t=1,2,\ldots,s.$ Thus $j_s$
exceeds $s$ elements of the $r-$th row. This forces $j_s>i_s,$ as required.
\QED
\par
If the entries 
of $S$ all exceed those of $T$ then
one only needs to shift numbered boxes (to the left). 
In a similar fashion $T\choose S$ is defined. One has 
$\sigma{ T\choose S}=(\sigma(T)\sigma(S)).$
\par
Note that $T^{i,i+1}={T^i\choose T^{i+1}},$ for example.
\parno
{\bf 2.8}\ \ \ 
Take $j \in {\Na}$ such that 
$j\not\in <R>.$ Let $a_i$
be the smallest entry of $R$ greater than $j$ and 
set:
$$(R+j):=(a_{\sr 1},\cdots,a_{i-\sr 1},j,a_i,\cdots)$$
\par
Define $w_r,\ w_c\in W$ by word forms
$$w_r(T):=[T^{\l_1^*},\cdots,T^1],\qquad
w_c(T):=[\hat T_1,\cdots, \hat T_{\l_1}]$$
\proclaim Lemma. {\rm [M5 or M1, 3.2.2]} One has 
$$Q(w_r(T))=Q(w_c(T))=T.$$
\par
\Pf
We give a proof for completion and to clarify the nature
of Robinson-Schensted correspondence.
\par
Let $T$
be a Young tableau and $b_1$ an integer $>0$ not belonging to
the entries of $T.$ The Robinson-Schensted (RS)
insertion of $b_1$ into $T$ gives a new Young tableau
$(T\da b_1).$ This is an inductive procedure in which $b_1,b_2,\ldots,$
are defined and $b_i$ is inserted into the $i-$th row of $T$.
In detail, write $T^1=(a_1,a_2,\ldots,a_l).$
If $b_1>a_l$ then $b_1$ is inserted into $T^1$ and the process stops,
that is 
$$(T\da b_1)=\pmat{ T^1+b_1\cr
                    T^{2,k}\cr}.$$
Otherwise let $i$ be the smallest integer such that $b_1<a_i.$
Then $(T\da b_1)^1$ is obtained by replacing $a_i$ by $b_1.$ Set $b_2=a_i.$ 
Then
$$(T\da b_1)=\pmat{ (T^1\da b_1)^1\cr
                    (T^{2,k}\da b_2)\cr}.$$
\par
Let $w=[a_1,a_2,\ldots,a_n]$ be a word form.
Then $Q(w)$ is obtained inductively as follows. Take $Q_0(w)$ to be the empty
tableau and set $Q_i(w)=(Q_{i-1}(w)\da a_i),\ \forall i:\ 1\leq i\leq n.$
Then $Q(w)=Q_n(w).$
\par
We prove lemma by induction on rows.
Let $T$ be a Young tableau and set $T^1=(a_1,a_2,\ldots,a_l).$
Then $a_i<a_{i+1},\ a_i<T^2_i,\ \forall i:\ 1\leq i\leq l$ (using the convention that the undefined
entries $T_i^j:\ j>2,$ of $T$ are set equal to $\infty$).
Conversely starting from $T^{2,k}$ and the sequence $R=(a_1,a_2,\ldots,a_l)$
with the above properties we obtain a Young tableau $S$ with $S^1=R,\ S^i=T^i,\ \forall i>1.$ 
One checks that $Q([w_r(T^{2,k}),R])=S$ and this proves the first part of the lemma.
The second obtains via $\sigma.$  
\QED
\parno
{\bf 2.9}\ \ \ 
Let us describe the jeu de taquin (cf. [Sch]) which removes $T^i_j$
from $T.$ The resulting tableau is denoted by $T-T^i_j$ and is 
obtained as follows.
\par
Consider a row $R=(a_{\sr 1},\cdots).$
Take 
$j \in {\Na}$ such that $j>a_{\sr 1}$ and 
$j\not\in <R>.$ Let $a_i$
be the greatest entry of $R$ smaller than $j$ and 
set:
$$(R \ua j):=(a_{\sr 1},\cdots,a_{i-1},j,a_{i+1},
                     \cdots).$$
One may remark that $((R \ua j)\da a_i)=R.$ 
\par 
Take  some element $a_i$ of $R$ and set:
$$(R-a_i):=(a_{\sr 1},\cdots,a_{i-\sr 1},a_{i+\sr 1},\cdots)$$
Similar operations are defined for a column $C$ except that 
we write $(C\la j)$ for 
\break
$\sigma(\sigma(C)\ua j).$
\par
With these preliminaries we write
\item{(1)} If $h(T^i_j)=1,$ then 
$$(T-T^i_j)=\pmat{T^{1,i-1}\cr
                          (T^i-T^i_j)\cr
                             T^{i+1,m}\cr}$$
\item{(2)} If $h(T^i_j)>1,$ then
\itemitem{(i)} If $T^i_{j+1}> T^{i+1}_j$ or $\l_i=j,$ set
$$(T-T^i_j)=\pmat{T^{1,i-1}\cr
                            (T^i\ua T^{i+1}_j)\cr
                            (T^{i+1,m}-T^{i+1}_j)\cr}$$
\itemitem{(ii)} If $T^i_{j+1}< T^{i+1}_j$ or $\l_j^*=i,$ set
$$(T-T^i_j)=(T_{1,j-1},(T_j\la T^i_{j+1}),
(T_{j+1,\l_{\sr 1}}-T^i_{j+1}))$$
\proclaim Theorem. {\rm [Sch]} If $T$ is a Young tableau then 
so is $(T-T^i_j).$
\par
\parno
{\bf 2.10}\ \ \ 
Given $w\in W,\ \ T\in\bT_n,\ \ \Vscr\in\bV,$ 
a standard parabolic subgroup $\bP$ and a standard parabolic
subalgebra $\gp.$ Let $\gp_{\a_i}$ be the minimal
standard parabolic subalgebra including $X_{-\a_i}$ and let $\bP_{\a_i}$
be the corresponding standard parabolic subgroup.
Define their $\tau$-invariants to be
$$ \eqalign{\tau(w)&    :=\Pi\cap S(w^{-1}),\cr
            \tau(T)&    :=\{\a_i\ : r_{\sr T}(i+1)>r_{\sr T}(i)\},\cr
            \tau(\bP)& :=\{\a_i\ :\bP_{\a_i}\st \bP\},\cr
            \tau(\gp)& :=\{\a_i\ :\gp_{\a_i}\st \gp\},\cr
            \tau(\Vscr)&:=\{\a_i\ : \bP_{\a_i}(\Vscr)=\Vscr\}.\cr}$$
Note that $\bP$ (resp. $\gp$) is uniquely determined by its 
$\tau$-invariant. Let $\bP_T=\bP_{\Vscr_T}$ be the stabilizer 
of $\Vscr_T.$ It is a standard parabolic subgroup of $\bG$ and
we set $\gp_{\sr T}=\gp_{\sr \Vscr_T}=\Lie P_T.$ By say [J,\ \S 9] one has
\proclaim Lemma. $\tau(w)=\tau(\Vscr_w)=\tau(Q(w))=\tau(\bP_{\Vscr_w})=
\tau(\gp_{\Vscr_w}).$
\par
\parno
{\bf 2.11}\ \ \ 
Given $\Iscr\st \Pi,$ let $\bP_{\Iscr}$ denote the unique 
standard parabolic subgroup of $\bG$ such that 
$\tau(\bP_{\Iscr})=\Iscr.$ Let $\bM_{\Iscr}$ be the unipotent 
radical of $\bP_{\Iscr}$ and $\bL_{\Iscr}$ a Levi factor.
Let $\gp_{\sr \Iscr},\ \gm_{\sr \Iscr},\ \gl_{\sr \Iscr}$ denote the 
corresponding Lie algebras. These notation will be conserved 
throughout though the subscripts may sometimes be dropped.
\par
Given $w\in W,$ let $\ell(w)$ denote its reduced length.
Let $W_{\Iscr}$ be the subgroup of $W$ generated by the 
$s_{\a}\ :\ \a\in \Iscr$  and $w_{\sr \Iscr}$ its unique longest element.
It is also defined by condition $S(w_{\sr \Iscr})=R^+\cap \Na \Iscr.$
The following is well known and easy to check.
\proclaim Lemma. $\Vscr_{\Iscr}:=\Vscr_{w_{\sr \Iscr}}$ is the unique orbital
variety with closure $\gm_{\sr \Iscr}.$ Moreover if $\tau(\Vscr)\supset \Iscr,$
then $\Vscr\st\ov{\Vscr_{\Iscr}}=\gm_{\sr \Iscr}.$
\par
\parno
{\bf Remark.}  Set $\Oscr_\Iscr=\bG\Vscr_\Iscr.$ One calls $\Vscr_\Iscr$ the Richardson
component (of $\Oscr_\Iscr\cap\gn$) defined by $\Iscr.$
With respect to the order relation defined in 2.6,
$\Vscr_{\Iscr}$ is the unique minimal orbital variety with $\tau-$invariant
$\Iscr.$ In general if $\gn\cap^y\gn\st \gn\cap^w\gn$ then
trivially $\Vscr_y\st \ov{\bBo(\gn\cap^w\gn)}=\ov \Vscr_w.$
By the lemma the converse holds if $w=w_{\sr \Iscr},$ for some $\Iscr\st\Pi;$
but fails in general [M5, 4.1.1].
\parno
{\bf 2.12}\ \ \  Given $\Iscr\st\Pi.$ Write $T_{\Vscr_\Iscr}$ simply as 
$T_\Iscr.$ It  
is obtained by the following rules
\item {(i)} $r_{\sr T_{\Iscr}}(1)=1.$
\item {(ii)} Given $s:\ 1<s\leq n$  
 $$r_{\sr T_{\Iscr}}(s)=\cases{1, &if $\a_{s-1}\not\in\Iscr,$\cr
                           r_{\sr T_{\Iscr}}(s-1)+1,& otherwise.\cr}$$
\par
A useful way to present $T_{\Iscr}$ is as follows. Partition
$\{1,2,\ldots,n\}$ into connected subsets 
$C_j:=\{b_j,b_j+1,\ldots,b_{j+1}-1\}$ by choosing a strictly 
increasing sequence 
\break
$1=b_1<b_2<\cdots<b_{l+1}=n+1.$
Setting $\Iscr=\{\a_i\ |\ i,i+1\ {\rm belong\ to\ some\ }C_j\}$
defines a bijection between the set of all such  partitions and 
the set of subsets of $\Pi.$ Given $\Iscr\st \Pi,$ let
$\{C_i^{\Iscr}\ :\ i=1,2,\ldots,l\}$ be the corresponding 
connected subsets which we view as columns. (Sometimes we may
omit the $\Iscr$ superscript.) Then in the notation of 2.7 we 
have $T_{\Iscr}=(C^{\Iscr}_1,C^{\Iscr}_2,\ldots,C^{\Iscr}_l).$
Of course this involves some sliding of boxes to the left.
However there are some advantages in this presentation. For example
$T_\Iscr -1$ is obtained by simply replacing $C_1$ by $C_1-1$ and
$T_\Iscr-n$ is obtained by simply replacing $C_l$ by $C_l-n.$ Again
one easily checks the
\proclaim Proposition. For all $\Iscr\st \Pi,$ the word form
of $w_{\sr \Iscr}$ is given by 
$$w_{\sr \Iscr}=[\widehat C^{\Iscr}_1,\widehat C^{\Iscr}_2,\ldots,
\widehat C^{\Iscr}_l].$$
\par
We call $C^\Iscr_i$ the $i-$th chain of $T_\Iscr$ and 
$T_\Iscr=(C^{\Iscr}_1, C^{\Iscr}_2,\ldots, C^{\Iscr}_l)$ the chain form of
$T_\Iscr.$
\parno
{\bf Example}
\parno
Consider $\Iscr=\{\a_1,\a_4,\a_5,\a_6,\a_9\}$ in $\gs\gl_{10}.$ Then
$$T_\Iscr =
\vcenter{
\halign{& \hfill#\hfill
\tabskip4pt\cr
\multispan{11}{\hrulefill}\cr
\ss
\vb & 1 &  &  3 & &  4 && 8 && 9 &\ts\vb\cr
\vs
&&&&&&\multispan{5}{\hrulefill}\cr
\ss
\vb & 2 &  &  5 && 10 & \ts\vb\cr
\vs
&&\multispan{5}{\hrulefill}\cr
\ss
\vb & 6 & \ts\vb\cr
\vs
&&\cr
\ss
\vb & 7 & \ts\vb\cr
\vs
\multispan{3}{\hrulefill}\cr}}$$
or in chain form $T_\Iscr=(C_1,C_2,C_3,C_4,C_5)$ where $C_1=\{1,2\},\ C_2=\{3\},\ 
C_3=\{4,5,6,7\},\ C_4=\{8\}$ and $C_5=\{9,10\}.$ 
\parno
{\bf 2.13}\ \ \
Given an orbital variety $\Vscr,$ set $\Oscr=\bG\Vscr.$ It is a 
nilpotent orbit. Given a nilpotent orbit $\Oscr\pr\st\ov \Oscr$
one may try to describe $\ov\Vscr\cap\Oscr\pr.$ In general it is 
not even known if this intersection is equidimensional. However if
$\Vscr=\Vscr_{\Iscr}$ for some $\Iscr\st\Pi$ a  complete 
answer given by [S1], [M5 or M1] and 2.11. Set 
$\Oscr_{\Iscr}=\bG\Vscr_{\Iscr}.$
\proclaim Theorem. Let $\Oscr$ be a nilpotent orbit in 
$\ov\Oscr_{\Iscr}=\bG\gm_{\sr \Iscr}.$ Then $\gm_{\sr \Iscr}\cap \Oscr$
is the union of orbital varieties $\Vscr$ in $\Oscr\cap\gn$
satisfying $\tau(\Vscr)\supset\Iscr.$
\par
\Pf
By [S1, p. 456, last corollary] $\gm_{\sr \Iscr}\cap\Oscr$
is equidimensional. 
\par
By [M1, 4.1.8 or M5, 3.5.6] it contains at least one orbital
variety and so $\dim \Vscr={1\over 2}\, \dim \Oscr,$ for every irreducible 
component $\Vscr$ of $\Oscr\cap\gm_{\sr \Iscr}.$ Yet $\Vscr\st \Oscr\cap\gn,$
which also has dimension ${1\over 2}\, \dim \Oscr.$ Hence $\Vscr$ is an orbital
variety associated to $\Oscr.$ Then by 2.11 every orbital variety $\Vscr$
associated to $\Oscr$ lies in $\gm_{\sr \Iscr}$ if and only if $\tau(\Vscr)\supset\Iscr.$
\QED
\parno
{\bf Remark.} By 2.10, $\a_i\in\tau(\Vscr_T)$
if and only if $r_{\sr T}(i+1)>r_{\sr T}(i).$ In particular if $\Vscr>\Vscr_{\Iscr}$
then for any $i$ one has $r_{\sr T_{\Vscr}}(i)\geq r_{\sr T_{\Iscr}}(i).$
\parno 
{\bf 2.14}\ \ \
Take $\Iscr\st\Pi$ and recall the notation of 2.11.
Set $\bBo_{\Iscr}=\bBo\cap\bL_{\Iscr},\ \gn_{\sr \Iscr}=\gn\cap\gl_{\sr \Iscr}$ and 
$W^{\Iscr}=\{w\in W\ |\ w\Iscr\st R^+ \}.$ 
We have decompositions $\bBo=\bM_{\Iscr}\spro\bBo_\Iscr,\ \gn=\gn_{\sr \Iscr}\oplus\gm_{\sr \Iscr},\
W^\Iscr\times W_\Iscr \raeq W.$ They define projections 
$\bBo\ra \bBo_\Iscr,\ \gn\ra\gn_{\sr \Iscr}$ and $W\ra W_\Iscr$
which we will denote by $\pi_{\sr \Iscr},$ or simply by $\pi.$ 
Let $\Vscr^\Iscr_{\pi(w)}$ be the orbital variety of $\gl_{\sr \Iscr}$
with closure $\ov{\bBo_\Iscr(\gn_{\sr \Iscr}\cap^{\pi(w)}\gn_{\sr \Iscr})}.$
\proclaim Proposition. {\rm [M1, 4.1.2 or M5, 4.2.2]} For all 
$w\in W,\ \Iscr\st \Pi$ one has $\pi_{\sr \Iscr}(\ov\Vscr_w)=
\ov{\pi(\Vscr_w)}=\ov{\Vscr^\Iscr_{\pi(w)}}.$
\par
\parno
{\bf Remark 1.} It is clear that $\pi_{\sr \Iscr}$ is inclusion 
preserving. Thus if $\Vscr_w,\ \Vscr_y$ are orbital varieties of
$\gg$ with $\ov \Vscr_w\st\ov\Vscr_y,$ then 
$\ov\Vscr_{\pi(w)}^\Iscr\st\ov\Vscr_{\pi(y)}^\Iscr.$
\parno
{\bf Remark 2.} Suppose $\ga\st\gn$ is $\gh$ stable. Since root 
subspaces are one-dimensional we obtain 
$\pi_{\sr \Iscr}(\ga)=\ga\cap\gn_{\sr \Iscr}.$ Again $\bM_\Iscr$ acts by 1 on 
$\gn/\gm_{\sr \Iscr}.$ Then since $\pi_{\sr \Iscr}$ is continuous
$$\pi_{\sr \Iscr}(\ov{\bBo\ga})\st
\ov{\pi_{\sr \Iscr}(\bBo\ga)}\st\ov{\bBo\ga\cap\gn_{\sr \Iscr}}\st
\ov{\pi_{\sr \Iscr}(\bBo\ga)}$$
and so $\ov\Vscr^\Iscr_{\pi(w)}=\ov\Vscr_w\cap\gn_{\sr \Iscr}.$
\par
For $1\leq i< j\leq n$ set $[i,j]=\{r\in{\Na}\ | i\leq r\leq j\},\ 
\Pi_{i,j}=\{\a_r\ |\ r,r+1\in [i,j]\}$ and $\pi_{i,j}=
\pi_{\sr \Pi_{i,j}}.$
Define $\bT_n\ra \bT_{j-i+1},$ through the jeu de taquin applied
to the entries (taken in any order) of $T$ not lying in $[i,j].$
By [M1, 2.4.16, 4.1.1 or M5 4.3.3] one has
$$\pi_{i,j}(\ov\Vscr_T)=\ov\Vscr_{\pi_{i,j}(T)}. \eqno {(*)} $$
\par
\proclaim {\ \ \ \ 2.15\ \ \ Definition}. Let $\geq$ be an order relation on $\bV.$
Given $\Vscr_{T_2}>\Vscr_{T_1}$
We call $\Vscr_{T_2}$ a descendant of $\Vscr_{T_1}$
(respectively $T_2$ a descendant of $T_1$) if for any
orbital variety $\Wscr$ such that 
$\Vscr_{T_2}\geq\Wscr\geq \Vscr_{T_1}$ one has 
$\Wscr=\Vscr_{T_2}$ or $\Wscr=\Vscr_{T_1}.$ One calls $\Vscr_{T_2}$ a geometric 
descendant of $\Vscr_{T_1}$ when $\geq$ is the geometric order. 
\par
If unqualified descendant will mean geometric descendant.
\par
The set of descendants of any $T_{\Iscr}\ :\ \Iscr\st\Pi$ is described explicitly in [M2, 2.6]
 The simplification that results in our case can be
understood as follows.
\par
Define the Duflo (or weak Bruhat) order on $W$ by $y\supseteq w$ if $\gn\cap^y \gn\st\gn\cap^w\gn.$
Through the map $w\mapsto Q(w)$ it induces an order relation $\supset$ on $\bT_n$ and hence on
$\bV$ called the (induced) Duflo order. Remark 2.11 states just that $\Vscr_y\st\ov\Vscr_w$ if $y\supseteq w,$
but the converse can fail unless $w=w_{\sr \Iscr}.$ The description
of the induced Duflo order on $\bT_n$ is itself a quite non-trivial problem; but the complete solution
was given in [M1, 3.3.3, 3.4.5 or M5, 3.3.3, 3.4.4]. By the above remark the set of geometric descendants of a Richardson
component is contained in the set of Duflo descendants in turn described by [M2, 2.6].
\par
Present $T_\Iscr$ as in 2.12, that is we write 
$T_{\Iscr}=(C^\Iscr_1,C^\Iscr_2,\ldots,C^\Iscr_l).$
Set $c_i=|C_i|.$ Then $\s_i=C^{c_i}_i$ is the largest entry of 
$C_i$ and moreover belongs to $T_\Iscr^{c_i}.$
\par
For any $1\leq i\leq l,$ suppose that there exists $j\ :\ 1\leq j<i$
such that $c_j\geq c_i.$ Take the maximal integer $i\pr\ :\ 1\leq i\pr<i$ such that 
$c_{i\pr}$ is minimal with the property that $c_{i\pr}\geq c_i.$
Define 
$$T_\Iscr(i):=\pmat{T^{1,c_i-1}\cr
                    T^{c_i}-\s_i\cr
                    T^{c_i+1,c_{i\pr}}\cr
                    T^{c_{i\pr}+1}+\s_i \cr
                    T^{c_{i\pr}+2,k} \cr}$$
In other words we move $\s_i$ from the $c_i-$th to the $(c_{i\pr}+1)-$th row according
to the rules of 2.9.
In example 2.12 one has $c_{\sr 2}=1>c_{\sr 1}=2;\ c_{\sr 3}=4$ is maximal, $c_{\sr 4}=c_{\sr 2}=1$ and $c_{\sr 5}=c_{\sr 1}=2.$
Thus our procedure defines 3 tableaux:
$$T_\Iscr(2) =
\vcenter{
\halign{& \hfill#\hfill
\tabskip4pt\cr
\multispan{9}{\hrulefill}\cr
\ss
\vb & 1 &  & 4  & & 8 && 9 &\ts\vb\cr
\vs
&&&&&&\multispan{3}{\hrulefill}\cr
\ss
\vb & 2 &  &  5 && 10 & \ts\vb\cr
\vs
&&&&\multispan{3}{\hrulefill}\cr
\ss
\vb & 3 && 6& \ts\vb\cr
\vs
&&\multispan{3}{\hrulefill}\cr
\ss
\vb & 7 & \ts\vb\cr
\vs
\multispan{3}{\hrulefill}\cr}},\ 
T_\Iscr(4) =
\vcenter{
\halign{& \hfill#\hfill
\tabskip4pt\cr
\multispan{9}{\hrulefill}\cr
\ss
\vb & 1 &  & 3  & & 4 && 9 &\ts\vb\cr
\vs
&&&&&&&&\cr
\ss
\vb & 2 &  &  5 && 8&& 10 & \ts\vb\cr
\vs
&&\multispan{7}{\hrulefill}\cr
\ss
\vb & 6& \ts\vb\cr
\vs
&&\cr
\ss
\ss
\vb & 7 & \ts\vb\cr
\vs
\multispan{3}{\hrulefill}\cr}},\
T_\Iscr(5) =
\vcenter{
\halign{& \hfill#\hfill
\tabskip4pt\cr
\multispan{11}{\hrulefill}\cr
\ss
\vb & 1 &  &  3 & &  4 && 8 && 9 &\ts\vb\cr
\vs
&&&&\multispan{7}{\hrulefill}\cr
\ss
\vb & 2 &  &  5 &\ts\vb\cr
\vs
&&&&\cr
\ss
\vb & 6 && 10 & \ts\vb\cr
\vs
&&\multispan{3}{\hrulefill}\cr
\ss
\vb & 7 & \ts\vb\cr
\vs
\multispan{3}{\hrulefill}\cr}}.
$$
\proclaim Lemma. $T_\Iscr(i)$ is a Young tableau.
\par
\Pf
Indeed the largest element $\s_i$ of the $i-$th chain is relocated on the $i\pr-$th chain where it becomes the largest element since $i\pr<i.$
Then the assertion follows from lemma 2.7.
\QED
\parno
{\bf 2.16}\ \ \ 
As noted in 2.12 the tableau $\pi_{1,n-1}(T_\Iscr)$
is obtained from $T_\Iscr$ by eliminating $n$ from $C_l.$ It corresponds 
to the subset $\Iscr_n=\Iscr\cap \Pi_{1,n-1}$
with $C^{\Iscr_n}_i=C^\Iscr_i\ :\ i<l$ and $C^{\Iscr_n}_l=
C^\Iscr_l-n.$ Similarly $\pi_{2,n}(T_{\Iscr})$ is obtained from
$T_\Iscr$ by eliminating 1 from $C_1$ and corresponds to 
the  subset $\Iscr_1=\Iscr\cap\Pi_{2,n}$ with 
$C^{\Iscr_1}_i=C^\Iscr_i\ :\ i>1$ and $C^{\Iscr_1}_1=C^\Iscr_1-1.$
\proclaim Lemma. Assume that $T_{\Iscr}(i)$ is defined. Then
$i>1$ and
\item {(i)} $\pi_{1,n-1}(T_\Iscr(i))=T_{\Iscr_n}(i),$ if $i<l.$
\item{(ii)} $\pi_{2,n}(T_\Iscr(i))=T_{\Iscr_1}(i),$ if $c_{\sr 1}\ne c_i$
or there exists $j\ :\ 1<j<i$ such that $c_j= c_i.$
\item{(iii)} $\pi_{2,n}(T_\Iscr(i))=T_{\Iscr_1},$ if 
$c_{\sr 1}=c_i$ and $c_j\ne c_i,\ \forall j:\ 1<j<i.$

\par
\Pf
(i) is an immediate consequence of $\s_i<n,$ for $i<l.$ 
The hypothesis of (ii) implies
that $T_{\Iscr_1}(i)$ is defined. If in the choice of $i\pr$ (as above) 
one can choose 
$i\pr\ :\ 1<i\pr<i,$
then both $T_{\Iscr}(i)$ and $T_{\Iscr_1}(i)$ are obtained 
by moving $\s_i$
from  $c_i-$th row to $(c_{i\pr}+1)-$th row, an operation which commutes 
with replacing $C_1$
by $C_1-1.$ If necessarily $i\pr=1$ but $c_{\sr 1}>c_i$ then $\s_i$ goes from 
$c_i-$th row to  
$(c_{\sr 1}+1)-$th row in forming $T_\Iscr(i),$ but then goes to $c_{\sr 1}-$th row in
forming $\pi_{2,n}(T_\Iscr(i)).$ However if $c_{\sr 1}>c_i$ this is exactly what happens in
forming $T_{\Iscr_1}(i)$ from $T_{\Iscr_1}.$ Hence (ii). Under the hypothesis of (iii),
$\s_i$ is pushed back into $c_i-$th row in forming 
$\pi_{2,n}(T_\Iscr(i))$ which is hence $T_{\Iscr_1}.$ 
\QED
\parno
{\bf 2.17}\ \ \ Retain the conventions of 2.15.
\proclaim Proposition. Take $\Iscr\st \Pi.$ Assume $T_\Iscr(i)$ is defined 
and hence $i\pr$ is defined. If $c_{i\pr}=c_i$ then $T_\Iscr(i)$ is a 
descendant of $T_\Iscr$ of codimension 1
and all the descendants of codimension 1 
of $T_\Iscr$ are so obtained.
\par
\Pf
By 2.15, $T=T_\Iscr(i)$ is a Young tableau and 
$\Vscr_T\st\ov\Vscr_\Iscr.$ 
\par
By 2.4 (2) we have 
$$\dim\Vscr_T=\dim\Vscr_\Iscr-|c_{i\pr}-c_i|+1. \eqno{(*)}$$
Take $c_{i\pr}=c_i,$ then $\Vscr_T$ is a descendant just by 
dimension considerations.
\par
Conversely assume that $T>T_\Iscr$ with $\Vscr_T$ of codimension one in $\gm_{\sr \Iscr}.$
It remains to show that $T=T_\Iscr(i),$ for some $i.$
\par 
By Remark 2.13 we have
$$r_{\sr T}(s)\geq r_{\sr T_\Iscr}(s),\qquad {\rm for\quad all\quad} s.$$
On the other hand by 2.4(2) and 2.5, $\sh T$ is obtained from
$\sh T_\Iscr$ by lowering exactly one box by row. Combined with the previous result this means 
that there is exactly one entry, say $m$ of $T_\Iscr$ which is displaced on passing to
$T$ and moreover it appears exactly one row below. Suppose $\a_m\in\tau(T_\Iscr).$
Since $\tau(T)\supset\tau(T_\Iscr)$ by 2.11, it follows by Remark 2.13, that
$$r_{\sr T}(m+1)\geq r_{\sr T}(m)+1=r_{\sr T_\Iscr}(m)+2,\quad {\rm whilst}\quad r_{\sr T_\Iscr}(m+1)=r_{\sr T_\Iscr}(m)+1.$$
This contradicts $r_{\sr T}(m+1)=r_{\sr T_\Iscr}(m+1).$ Thus, in the previous notation $m=\s_i,$ 
the latter being the longest
element of some chain $C_i$ of $T_\Iscr$ of length $c_i.$ By the  above we  also have 
$T^{c_i}=T_\Iscr^{c_i}-m$ and $T^{c_i+1}=T_\Iscr^{c_i+1}+m$ whilst the remaining rows of $T,T_\Iscr$ 
must coincide. Finally we must have some $i\pr<i$ with $c_{i\pr}=c_i,$ for otherwise $m$
would be strictly less then an entry in $T$ above it. We conclude that $T=T_\Iscr(i),$
as required. 
\QED
\parno 
{\bf 2.18}\ \ \ By 2.17 we conclude that the orbital varieties of codimension 1
in $\Vscr_\Iscr$ are obtained by moving $\s_i$ from row $c_i$ to row $(c_i+1)$ wherever
$|T^{c_i}_\Iscr|\geq |T^{c_i+1}_\Iscr|+2.$ Indeed if $\s_i$ belongs to $T_j$
(and here $j\leq i$ because some sliding of boxes from the chain form of $T_\Iscr$ may be 
necessary) then the 
condition $c_{i\pr}=c_i$ implies that $|T_{j-1}|=|T_j|.$
\par  
More precisely there are exactly $t+1:=\sharp \{j\ |\ 1\leq j\leq i\ {\rm with}\ c_j=c_i\}$ columns of $T$,
namely $T_j,T_{j-1},\ldots,T_{j-t},$ having length $|T_j|.$ In addition suppose
that $i$ is maximal for a given value $c_i$ that is $c_j\ne c_i$ for any $j>i.$
Then there are exactly $t$ choices of $s$ such that $T_\Iscr(s)$ is defined and has the same shape
as $T_\Iscr(i).$
\par
Observe that the partition $\l_{\sr \Iscr}=\sh T_\Iscr$ is defined by chain lengths namely we
have $\l^*_{\sr \Iscr}=\{c_i\ :\ i=1,2,\ldots,l\}$ appropriately ordered. Fix $c\in\l^*_{\sr \Iscr}$
and set $t+1=\sharp \{i\ |\ c_i=c\}.$ Suppose $t\geq 1$ and let $\Oscr_\Iscr(c)$ be defined through
the partition $\l_{\sr \Iscr}(c)$ obtained from $\l_{\sr \Iscr}$ by replacing the appropriate pair $(c,c)$ in $\l_{\sr \Iscr}^*$
by $(c+1,c-1).$ Then we may summarize the above by

\proclaim Proposition. Fix $\Iscr\st\Pi$ and $c\in{\Na}^+.$
Set $t=\sharp \{i\ |\ c_i=c\}-1.$ Assume $t\geq 1.$ Then $\gm_{\sr \Iscr}\cap\Oscr_\Iscr(c)$ is a union
of $t$ orbital varieties of codimension 1 in $\gm_{\sr \Iscr}.$ Moreover these orbital varieties correspond to 
the standard tableaux $T_\Iscr(i)\ :\ c_i=c,$ except for the minimal $i$ that is $i:\ c_i=c$ and $c_j\ne c$
for all $j<i.$  
\par
Note that orbital varieties of codimension 1 in the nilpotent radical of a parabolic appear exactly when 
there are repetitions in the dual partition associated to this parabolic.
\parno
{\bf 2.19}\ \ \ Let $\gm^-$ be the opposed algebra of $\gm=\gm^+$ identified with $\gm^*$
through the Killing form. By Krull's theorem the ideal $I(\ov\Vscr)$ of definition of an orbital variety 
closure $\ov\Vscr$ in the conclusion of 2.18 is principal in $S(\gm^-)$ and we denote by
$f_{\sr \Iscr}(i),$ or simply $f,$ the corresponding irreducible generator in $S(\gm^-).$
\par
The form of $f_{\sr \Iscr}(i)$ was conjectured by Benlolo and Sanderson in [BS].
We refer to this and their conjecture simply as BS. We prove this conjecture in 
Sect. 3.
\par
As explained in 1.4, 1.5 the general form of $f_{\sr \Iscr}(i)$ is $\gr a$ for some minor $a$ of a generic
matrix in $\gm_{\sr \Iscr} +\Id.$ 
They took a suitable minor lying in a bottom left-hand corner. In this they
needed results similar to [M5] to reduce to the case when $\Iscr,\ c$ are chosen so that $c=c_{\sr 1}=c_l$
and $c_i\ne c$ for $1<i<l.$ Then their conjecture is equivalent to saying that $a$ is just the 
$(n-c)\times (n-c)$ minor lying in the bottom left-hand corner.
\par
The BS reduction can be read off from 2.16 and 2.18 which makes precise the result of [M5]
which one needs and we believe clarifies their analysis. A key point is the knowledge of
the orbital varieties of codimension 1 in $\gm_{\sr \Iscr}$ given by 2.17 developed from the results 
of [M5]. We remark that the codimension 1 case is relatively easy and barely uses the full power of
[M5].
\par
For all $1\leq i,j\leq n,$ set $X_{i,j}={\Co}e_{i,j}.$ Given $w\in W,$ we give the affine
\break 
$|R^+|-\ell(w)$ dimensional space $X(w)$ a matrix presentation through $X(w)_{i,j}=X_{i,j},$ wherever $i<j$
and $w(i)<w(j).$
\par
Let $x_{j,i}$ denote the coordinate function on $\gg$ defined by
$$x_{j,i}(e_{r,s})=\cases{-1, & if $(r,s)=(i,j),$\cr
                          0, & otherwise.\cr}$$
Then the Poisson bracket $\{,\}$ defined on $S(\gg^*)$ through the Lie bracket on $\gg$ satisfies
$$\{x_{i,j},x_{r,s}\}=\d_{j,r}x_{i,s}-\d_{s,i}x_{r,j} \eqno{(*)} $$
where $\d_{i,j}$ is the Kronecker delta. Setting $x_{i,j}=e_{i,j}$ for $i>j$ identifies $\gm^-$ with $\gm^*.$
\par
Given $w\in W,$ let $\M(w)$ be the matrix with entries
$$\M(w)_{i,j}=\cases{x_{i,j},& if $X_{j,i}\st X(w),$\cr
                     0, & otherwise.\cr}$$
When $w=w_{\sr \Iscr},$ we set $\M(w)=\M_\Iscr.$
\par
Assume $c_{\sr 1}=c_l=:c$ and let $M_\Iscr^c(t)$ denote the $(n-c)\times(n-c)$ minor in the bottom left-hand 
corner of $\M_\Iscr+t\Id,$ that is
$$M_\Iscr^c(t):= 
\left|\matrix{x_{c+1,1} & \cdots\cdots &x_{c+1,c}&t   &0      &\cdots&0\cr
                              x_{c+2,1} & \cdots\cdots &\ldots   &\ldots    &t&\cdots&0\cr
                              \vdots    & \vdots       &\vdots   &\vdots    &\vdots &\vdots&\vdots\cr
                              x_{n-c,1} & \cdots\cdots &\ldots   &\ldots    &\ldots &\ldots&t\cr
                              \vdots    & \vdots       &\vdots   &\vdots    &\vdots &\vdots&\vdots\cr
                              x_{n,1}   &\cdots\cdots  &x_{n,c}  &x_{n,c+1} &\ldots &\ldots&x_{n,n-c}\cr}\right|$$
Note that there are some zeros in the dot places of the determinant in correspondence with the definition
of matrix $\M_\Iscr.$
\parno
{\bf Example}
\parno
Consider $\Iscr=\{\a_{\sr 1},\ \a_{\sr 4}\}$ in $\gs\gl_{\sr 5}.$
In that case $c_{\sr 1}=c_{\sr 3}=2$ and 
$$\M_\Iscr=\left[ \matrix{0          &0          &0          &0          &0\cr
                          0          &0          &0          &0          &0\cr
                          x_{\sr 3,1}&x_{\sr 3,2}&0          &0          &0\cr
                          x_{\sr 4,1}&x_{\sr 4,2}&x_{\sr 4,3}&0          &0\cr
                          x_{\sr 5,1}&x_{\sr 5,2}&x_{\sr 5,3}&0          &0\cr}\right]\quad {\rm and}\quad
M_\Iscr^2(t):= 
\left|\matrix{x_{\sr 3,1}&x_{\sr 3,2}&t\cr
              x_{\sr 4,1}&x_{\sr 4,2}&x_{\sr 4,3}\cr
              x_{\sr 5,1}&x_{\sr 5,2}&x_{\sr 5,3}\cr}\right|$$ 
\par
Developing $M_\Iscr^c(t)$ in powers of $t$ one obtains
$$M_\Iscr^c(t)=m_{n-c}+m_{n-c-1}t+\cdots+m_ct^{n-2c}.$$

Set 
$$d_i=\cases{c_i-c & if $c_i>c$\cr
             0     & otherwise\cr},\qquad d_{\sr \Iscr}=\sum\limits_{i=0}^k d_i,\qquad l_{\sr \Iscr}=n-c-d_{\sr \Iscr}$$
By [BS, lemma 2] the lowest non-zero coefficient is $m_{l_{\sr \Iscr}}$ that is the coefficient
of $t^{d_{\sr \Iscr}}.$  
The BS conjecture states that $f_{\sr \Iscr}(l)=m_{l_{\sr \Iscr}}.$ 
\par
In our example $d_{\sr \Iscr}=0$ and $l_{\sr \Iscr}=3$ so that 
$$m_{l_{\sr \Iscr}}=x_{\sr 3,1}x_{\sr 4,2}x_{\sr 5,3}+
x_{\sr 3,2}x_{\sr 4,3}x_{\sr 5,1}-x_{\sr 3,1}x_{\sr 4,3}x_{\sr 5,2}-x_{\sr 3,2}x_{\sr 4,1}x_{\sr 5,3}.$$

It is clear that $m_{l_{\sr \Iscr}}=\gr M^c_\Iscr(1).$
We do not need to know the explicit power of $t$ which divides $M_\Iscr^c(t),$
though this is used in one of our three proofs of irreducibility.
The BS conjecture can be expressed as the following theorem which we prove in the next section.
\proclaim Theorem. Take $\Iscr\st \Pi$ and let $T_\Iscr=(C_1,C_2,\ldots,C_l)$ be its chain form.
Set $c_i=|C_i|,$ for all $i.$
\item{(i)} Assume $c_{\sr 1}=c_l=:c$ and 
$c_j\ne c,\ \forall j:\ 1<j<l.$ Then $I(\ov\Vscr_{T_\Iscr(l)})$
is generated by $\{x_{ij}:\ X_{ij}\in\gp_{\sr \Iscr}\cap\gn^-\}$
and $m_{l_{\sr \Iscr}}.$
\item{(ii)} Suppose there exist $i<j$ such that $c_i=c_j$ and 
$c_s\ne c_i,\ \forall s:\ i<s<j.$ Set $u=C_i^1,\ v=C_j^{c_i}$
and $\Iscr\pr=\Iscr\cap\Pi_{u,v}.$ Then $I(\ov\Vscr_{T_\Iscr(j)})$
is generated by
\break 
$\{x_{ij}:\ X_{ij}\in\gp_{\sr \Iscr}\cap\gn^-\}$ and $m_{l_{\Iscr\pr}}$
defined with respect to $\Vscr_{T_{\Iscr\pr}}$ viewed as an orbital variety in  $\pi_{u,v}(\gn).$
\par
Note that the part (i) of the theorem is a special case of part (ii) and part (ii)
will be deduced from part (i).
\parno
{\bf 2.20}\ \ \ Let $i$ be a positive integer $<{n\over 2}$ and set
$$\b_i=\sum_{j=i}^{n-i}\a_j$$
Take $c$ as in 2.19 and set
$$\g_c=\sum_{i=1}^c\b_i.$$
As noted in [BS, VIII] one may check that $M_\Iscr^c(t),$
and hence $m_{l_{\sr \Iscr}}$ has a weight $-\g_c.$ Moreover with respect to the Cartan
inner product $(\ ,\ )$ one has 
$(\g_c,\a)=0,$
\break 
$\forall \a\in\Pi\setminus \{\a_c,\a_{n-c}\}.$
In particular
$$ (\g_c,\a)=0,\qquad \forall \a\in\Iscr. \eqno{(*)}.$$
\bigskip
\parno
{\bf 3. Proof of the Benlolo-Sanderson conjecture.}
\parno
{\bf 3.1}\ \ \  Let us recall some general facts about ideals of definition of orbital 
variety closures.
\par
Retain the notation of 2.11 and 2.18. Let $\bP$ be some standard parabolic subgroup
of $\bG,\ \gp=\Lie(\bP)$ and  $\gm$ its nilradical. Let $\Vscr$ be a closed subvariety of $\gm$ 
and $I$ (resp. $J$) its ideal of definition in $S(\gm^-)$ (resp. $S(\gg)$). One has 
$J=I+S(\gg)\gp$ and $J\cap S(\gm^-)=I$ from which one easily checks the
\proclaim Lemma. The following are equivalent
\item {(i)} $\{I,I\}\st I$ and $\{\gp,I\}\st J.$
\item {(ii)} $\{I,I\}\st I$ and $\Vscr$ is $\bP$ invariant.
\item {(iii)} $\{J,J\}\st J.$
\par
\parno
{\bf 3.2}\ \ \ One calls $\Vscr$ involutive if any one of the above holds.
If $\Vscr$ is involutive then so are its irreducible components.
\par
Suppose $\Vscr\st \gm$ is irreducible. Then $\ov{\bG\Vscr}$ contains a unique 
dense orbit $\Oscr.$ If in addition $\Vscr$ is involutive then 
$\dim (\Vscr\cap \Oscr)  \geq {1\over 2} \dim\Oscr$ since $\Oscr$ is a symplectic variety.
Yet $\dim(\gn\cap \Oscr)={1\over 2}\dim \Oscr$ and so $\Vscr\cap\Oscr$ must be an irreducible 
component of $ \gn\cap \Oscr$ hence an orbital variety associated to $\Oscr$ with 
closure $\Vscr.$
\parno
{\bf 3.3}\ \ \ Let $J_o$ be an ideal of $S(\gg)$ with radical $J.$ Suppose $\{J_o,J_o\}\st J_o.$ 
It is generally false that this implies $\{J,J\}\st J,$ so one cannot conclude that the zero 
variety of $J_o$ is involutive. This is a delicate point in general. In the present situation,
this difficulty can be avoided. Call $f\in S(\gm^-)$ multilinear if its degree in any one of the 
variables $x_{i,j}$ is at most 1. This property passes to the irreducible factors $f_1,f_2,\ldots,f_t$ 
and we let $\Vscr(f_1),\ldots,\Vscr(f_t)$ denote their sets of zeros in $\gm.$
Consider the condition 
$$\{\gp,f\}\st S(\gm^-)f+S(\gg)\gp. \eqno{(*)}$$
\proclaim Proposition. Suppose $f\in S(\gm^-)$ is multilinear. Then 
\item {(i)}  $S(\gm^-)f$ is semiprime.
\item {(ii)} $S(\gm^-)f+S(\gg)\gp$ is semiprime.
\item {(iii)} If $(*)$ holds then the $\Vscr(f_i)$ are orbital varieties of 
codimension 1 in $\gm.$
\par
\Pf
(i) is immediate. (ii) is easily deduced from (i). Finally (iii) follows from (ii) and 3.1, 3.2.
\QED
\parno
{\bf 3.4}\ \ \ Our proof of the BS conjecture involves three steps. The first is to show
that $f=m_{l_{\sr \Iscr}}$ satisfies 3.3 $(*).$ This is rather obvious especially from the quantum viewpoint 
of [J2, Lecture 7]. Nevertheless a delicate point is that $M_\Iscr^c(t)$ itself does not satisfy $(*).$
We may view its zero set as a deformation of $\Vscr(f);$ but which is not itself orbital nor $\bP$ 
invariant. (Our original motivation for  such deformation came from trying to define an Enright 
functor on orbital varieties itself inspired by the algorithm 
for $\Ann_{U(\gn^-)}v$ in [J1, 8.4]
based on the Enright functor.) This result is given in 3.11.
\par
The second step is to show that $f$ is irreducible. This is a delicate point. However we can give 
three different proofs. The first, given in 3.15, is a fairly explicit but includes some computations involving
the precise knowledge of $M_\Iscr^c.$ The second, outlined in 4.3, uses representation theory. We construct a simple
highest weight module $L$ with integral highest weight which is a strong quantization of $\Vscr(f).$
This is of interest in its  own right. Then we use the linear independence of the characteristic 
polynomials of orbital varieties and the irreducibility of the associated variety of $L$ which holds
[M4] in type $A.$ The third method, given in 4.7, is the least computational; but the most sophisticated. Here we 
construct a strong quantization; but not necessary having integral highest weight. Then we use 
the difficult fact that $V(\Ann L)$ is the closure of a nilpotent orbit. In type $A$ this has the 
relatively easy proof using mainly that orbit closures are 
normal [BK]. Finally we apply 2.18. In all these methods
it is crucial to use that $c_i\ne c,\ \forall 1<i<l.$ Otherwise $\Vscr(f)$ has exactly 
$t:=\sharp\{i\ |\ c_i=c\}-1$ components. A computational proof in some special cases
when $m_{l_{\sr \Iscr}}=M^c_\Iscr(0)$ was also given in [BS]. Even that is not trivial though this case can be viewed as
an easy consequence of 2.4 (1) combined with 2.18.
\par
Step three is to show that $M_\Iscr^c$ vanishes on $X(w)$ for some $w\ :\ Q(w)=T_\Iscr(l).$ Then $\bP_\Iscr$
invariance implies $\Vscr(f)\sp \ov{\bBo(\gn\cap^w\gn)}.$ This is shown in 3.8. Since both have codimension 1 
in $\gm_{\sr \Iscr},$ irreducibility finishes the proof. Notice this does not use involutivity; but by 3.1 the latter is equivalent 
to $\bP_\Iscr$ invariance. One may also avoid this last step by combining 3.3. and 2.18; but this is 
less satisfying.
\parno
{\bf 3.5}\ \ \ Given $w\in \bS_n$ with word form $[a_1,a_2,\ldots, a_n],$ let $w-s\in \bS_{n-1}$
be defined by deleting $s$ and standardizing the word form. For any matrix $\M$ let $\M^{i,j}$
be the matrix obtained from $\M$ by deleting the $i-$th row and $j-$th column. Recall the definition of
$X(w)$ given in 2.19.
\proclaim Lemma. For all $w\in W,\ s\in\{1,2,\ldots,n\}$ one has $X(w-s)=X(w)^{s,s}.$
\par
\Pf
Take $1\leq i<j\leq n.$ Then $w(\a_{i,j})\in R^+$ if and only if  $w(i)<w(j)$ with a similar assertion for $y=w-s.$
Yet recalling that we are standardizing $w-s$ we have up to standardization 
$$y(r)=\cases{a_r, & if $r<w^{-1}(s),$\cr
              a_{r+1},& if $r\geq w^{-1}(s),$\cr}$$
From this the assertion readily follows.
\QED
\parno
{\bf 3.6}\ \ \ Suppose $T=Q(w).$ It is generally false that $Q(w-s)=T-s,$ though it
is true for the canonical elements $w_r(T)$ and $w_c(T)$ defined in 2.8.
Here we shall only need to describe $Q(w_r(T)-s),$ when $s\in T^k,$ that is
to say when $s$ lies in the last row of $T.$ Recall that 
$w_r(T)=[T^k,T^{k-1},\ldots,T^1].$
\proclaim Lemma. One has $w_r(T)-T^k_j=w_r(T-T^k_j)$ for every entry
$T^k_j$ of the last row $T^k$ of $T.$
\par
\Pf
Indeed
$$(T-T^k_j)=\pmat{T^{1,k-1}\cr
                          (T^k-T^k_j)\cr}$$
which gives the required assertion.
\QED
\parno
{\bf 3.7}\ \ \ Fix $\Iscr\st\Pi$ and $s\in\{1,2,\ldots,n\}.$
Define $\Iscr\pr\st\Pi\setminus\{\a_{s-1}\}$ by
\break
$\Iscr\pr=\{\a_i\ |\ i<s-1,\ \a_i\in \Iscr\}\cup
\{\a_{i-1}\ |\ i\geq s,\ \a_i\in \Iscr\}.$ Write 
$T_\Iscr=(C_1^\Iscr,C_2^\Iscr,\ldots,C^\Iscr_l)$ as in 2.12.
\proclaim Lemma. 
\item{(i)} $w_{\sr \Iscr}-s=w_{\sr \Iscr\pr}.$
\item{(ii)} $Q(w_{\sr \Iscr}-s)=T_{\Iscr\pr}.$
\item{(iii)} Suppose $s\in C_j^\Iscr.$ Then up to standardization
$T_{\Iscr\pr}=(C_1^{\Iscr\pr},C_2^{\Iscr\pr},\ldots,C_l^{\Iscr\pr})$
where 
$$C_i^{\Iscr\pr}=\cases{ C^\Iscr_i, & if $i\ne j,$\cr
                         C^\Iscr_j-s, & if $i=j.$\cr}$$
\item {(iv)} If $s$ lies in the last row of $T_\Iscr$ then $T_{\Iscr\pr}=(T_\Iscr-s)$
after standardization.
\par
\Pf
Recall from 2.15 the notation $c_i:=|C_i|$ and $\s_i:=C_i^{c_i}.$ One has
$$\Iscr=\Pi\setminus \bigcup_{i=1}^{l-1}\a_{\s_i}$$
Then (i) is clear if $s\ne \s_i$ for any $i.$ Otherwise $\a_s\not\in\Iscr$
but also $\a_{s-1}\not\in\Iscr\pr.$ Hence (i). (ii) follows from (i) by 2.12.
(iii) is obtained from 2.12 just by definition of chains. Finally (iii)
implies (iv) since in that case $T_{\Iscr\pr}$ is obtained by eliminating
$s$ and sliding the rest to the left.
\QED
\parno
{\bf 3.8}\ \ \ Retain the notation of 2.19 and assume $c=c_{\sr 1}=c_l.$ Then $T_\Iscr(l)$ 
is defined. Set $M=M_\Iscr^c(t).$
\proclaim Proposition. Set $w=w_r(T_\Iscr(l)).$ Then $M^c_\Iscr(t)$ vanishes on 
$\gn\cap^w\gn,$ that is $M(X(w))= 0.$
\par
\Pf
The proof is by induction on $d_{\sr \Iscr}$ defined in 2.19.
If $d_{\sr \Iscr}=0,$ then $c_i\leq c,\ \forall 1\leq i\leq l$ and so $n$ occurs in the last
row of $T_\Iscr.$ Moreover $T_\Iscr(l)$ has $c+1$ rows with $n$ the unique entry in row $c+1.$
Consequently the word form of $w$ starts with $n.$ By 2.3 this forces $X_{i,n}\not\st 
\gn\cap^w\gn,\ \forall i<n.$  Thus the entries $x_{n,i}\ :\ i=1,2,\ldots,n-c,$ of the last row
of $M$ vanish on $X(w)$ and consequently so does $M.$
\par
Let us assume the assertion holds for $d_{\sr \Iscr}=d$ and take $d_{\sr \Iscr}=d+1.$ In particular there is some
$c_i>c.$ Let $k$ be the number of rows of $T:=T_\Iscr.$ Since $c_i>c$ this is also the number 
of rows of $T\pr:=T_\Iscr(l)$ and so the word form of $w$ starts with $s=T^k_1,$ where $s\not\in <C_1>.$
By 2.3 we conclude that $X_{i,s}\not\st \gn\cap^w\gn,$ for $i<s$ Thus the entries $x_{s,i}\ :\ i<s$
in $M^c_\Iscr(t)$ vanish on $X(w).$ This leaves $t$ as the only non-zero entry of row $s.$
Consequently
$$M(X(w))=tM(X(w)^{s,s}) \eqno{(*)}$$
it being understood that $n$ is reduced by 1 in defining the right hand side.
\par
By 3.5 one has $X(w)^{s,s}=X(w-s)$ whilst by 3.6 one has $X(w-s)=X(w_r(T\pr-s)).$ Recall the definition
of $\Iscr\pr$ from 3.7.
By 3.7 (iv), $T-s=T_{\Iscr\pr}.$ Then by 2.17  $T\pr-s$ is a descendant of $T_{\Iscr\pr}.$  
Moreover $d_{\Iscr\pr}=d_{\sr \Iscr}-1=d$ and $C_1^{\Iscr\pr}=C_1^\Iscr,\ C_l^{\Iscr\pr}=C_l^\Iscr$
up to standardization. In particular $|C_1^{\Iscr\pr}|=|C_l^{\Iscr\pr}|=c.$ 
Set $w\pr=w_r(T\pr -s),\ M\pr=M_{\Iscr\pr}^c(t).$ By the above 
and the induction hypothesis 
$M\pr(X(w\pr))= 0.$ Then by  $(*)$ and the above $M(X(w))= 0.$
\QED
\parno
{\bf 3.9}\ \ \ Retain the notation and hypotheses of 3.8. For $1\leq i,j\leq n,$ let $M^{i,j}$
denote $ij-$th cofactor of $M.$ Set $d=d_{\sr \Iscr}$ and recall that $t^d$ divides $M.$
The following is clear by row expansion of the determinant.
\proclaim Lemma. Choose $i,j$ such that $M_{i,j}=x_{c+i,j}\in 
\gm^-_{\sr \Iscr}.$ Then $M^{i,j}$ is divisible by $t^d.$
\par
\parno
{\bf 3.10}\ \ \ Take $f\in S(\gm^-).$ In computing $\{x_{i,j},f\}$ via 2.19 $(*)$ we obtain two types 
of terms. The first (resp. second) is obtained from the first (resp. second factor on the right-hand side.
We call them the terms obtained from $j$ (resp. $i$).
\proclaim Lemma. Suppose $i\in \{1,2,\ldots,n-1\}\setminus \{c,n-c\}.$
Then 
$$\{x_{i,i+1},M\}\st S(\gm^-_{\sr \Iscr})M+S(\gg)\gp_{\sr \Iscr}.$$
\par
\Pf
Let $\Xscr$ be matrix with entries $\Xscr_{i,j}=x_{i,j}$ and $X$ 
the $(n-c)\times (n-c)$ minor 
in the bottom left-hand corner of $\Xscr+t\Id.$ Then $M-X\in S(\gg)\gp_{\sr \Iscr}.$ Since $\gp_{\sr \Iscr}$ is a 
subalgebra it is enough to prove the corresponding assertion for $X.$
\par
Consider $\{x_{i,i+1},X\}.$ 
\item {(i)} Suppose $i>c.$ Then the sum of terms coming from $(i+1)$ forms a determinant
$X^{i+1}$ with the same rows as $X$ except that the entries $x_{i+1,s}+\d_{i+1,s}t$ on the $(i+1-c)-$th
row of $X$ are replaced by $x_{i,s}.$ Hence this term equals $-tX^{i+1-c,i}.$ 
\item {(ii)} Suppose $i<n-c.$ 
Then the sum of terms coming from $i$ forms a determinant $X_i$ with the same columns as $X$ except that the 
entries $x_{s,i}+\d_{s,i}t$ on the $i-$th column of $X$ are replaced by $-x_{s,i+1}.$ Hence this term
equals $tX^{i+1-c,i}.$
\par
If $c<i<n-c$ then terms cancel. If $i<c$ or $i>n-c$, both terms are zero. This proves the lemma.
\QED
\parno
{\bf 3.11}\ \ \ Let $f=m_{l_{\sr \Iscr}}.$
\proclaim Lemma. One has $\{\gp_{\sr \Iscr},f\}\st S(\gm_{\sr \Iscr}^-)f+S(\gg)\gp_{\sr \Iscr}.$
\par
\Pf
By 2.20 (*) $f$ has weight zero with respect to the Cartan 
sublagebra $[\gl_{\sr \Iscr},\gl_{\sr \Iscr}]\cap\gh$
of the semisimple Lie algebra $[\gl_{\sr \Iscr},\gl_{\sr \Iscr}].$ Since $S(\gg)$ and $S(\gg)\gl_{\sr \Iscr}$ are locally finite
$\gl_{\sr \Iscr}$ modules it follows from the theory of finite dimensional $\gl_{\sr \Iscr}$ modules that it suffices 
to show $\{x_{i,i+1},f\}\st S(\gm^-_{\sr \Iscr})f+S(\gg)\gp_{\sr \Iscr},$ for all $1\leq i\leq n-1$ and
\break
$\{x_{i+1,i},f\}\st S(\gm^-_{\sr \Iscr})f+S(\gg)\gp_{\sr \Iscr},$ for all $\a_i \in \tau (\gp_{\sr \Iscr}).$
By 3.10 among all $\{x_{i,i+1},f\}$ it remains 
to consider the cases $i=c,n-c.$ Both are similar and we consider only the first. Consider
$\{x_{c,c+1},M\}.$ The sum of terms coming from $c+1$ all lie in $S(\gg)\gp.$
The sum of terms coming from $c$ equal 
to $tM^{1,c}\ {\rm mod}\ S(\gg)\gp.$ Yet 
$M_{1,c}=x_{c+1,c}\in \gm_{\sr \Iscr}^-$ so by 3.9 this expression 
is divisible by $t^{d+1}\ {\rm mod}\ S(\gg) \gp_{\sr \Iscr}$
Thus $\{x_{c,c+1},f\}=0\ {\rm mod}\ S(\gg) \gp_{\sr \Iscr}.$
\par
As for $\{x_{i+1,i},f\}$ we show  that $\{x_{i+1,i},M\}=0$ 
exactly in the same manner as in 3.10. Let us sketch the proof. 
The cases $i<c$ or $i+1>n-c$ are the same so let us show this
for $i+1>n-c.$ In this case the sum
of terms coming from $i+1$ is zero  and the sum
of terms coming from $i$ forms
a determinant $M^i$ with the same rows as $M$ except that the 
$(i-c)$-th row is replaced by the $(i-c+1)$-th row so that $M^i=0.$
Now if $c<i<n-c-1$ then for all $s\ne i,i+1$ one has $M_{i-c,s}=0$
iff $M_{i+1-c,s}=0.$ As well for all $s\ne i-c,i+1-c$ one has 
$M_{s,i}=0$ iff $M_{s,i+1}=0$ Using this we get exactly as in 3.10 that 
the sum of terms
coming from $i$ results in $-tM^{i-c,i+1}$ and the sum of terms coming 
from $(i+1)$ results in $tM^{i-c,i+1}$ so that the terms cancel.
\par
\QED
\parno
{\bf 3.12}\ \ \ 
We may summarize the consequence of 3.3 (ii), 
3.8, 3.11 as follows. Adopt the notation and hypotheses of 2.19 (i).
In particular set $T\pr=T_\Iscr(l).$
\proclaim Proposition. The zero variety $\Vscr(m_{l_{\sr \Iscr}})$ of\ \ $m_{l_{\sr \Iscr}}$\ \  in\ \ $\gm_{\sr \Iscr}^-$\ \ 
is a union of orbital variety closures of codimension 1 one of which is $\Vscr_{T\pr}.$
\par
\parno
{\bf Remark.} We have not yet used that $c_j\ne c$ for $1<j<l.$ This is needed for irreducibility of $\Vscr(m_{l_{\sr \Iscr}}).$
This would follow from 2.18 and the irreducibility of $\bG\Vscr(m_{l_{\sr \Iscr}});$ but we only have a direct proof of the 
latter using representation theory (cf. 4.7).
\parno
{\bf 3.13}\ \ \  Let us show that (ii) of Theorem 2.19 results from (i). This obtains from the following general
remark. Take $\Iscr\st \Pi$ and $T_\Iscr=(C_1,\ldots,C_l)$ be its chain form.
Fix $i,j\ :\ 1\leq i\leq j\leq l$ and set $T_{\Iscr\pr}=(C_i,\ldots,C_j).$ Let $T\pr\in \bT_{|T_{\Iscr\pr}|}$
satisfy $T\pr>T_{\Iscr\pr}$ and set $T=(C_1,\ldots,C_{i-1},T\pr,C_{j+1},\ldots,C_l).$ Then $T>T_{\Iscr}$
and $\pi(T_\Iscr)=T_{\Iscr\pr},\ \pi(T)=T\pr,$ where $\pi=\pi_{\s_{i-1}+1,\s_j}.$ Again by 2.14 $(*)\ \ 
\pi(\Vscr_T)=\Vscr_{T\pr}$ and by 2.4 (2) the codimension of $\Vscr_T$ in $\gm_{\sr \Iscr}$ equals the codimension of 
$\Vscr_{T\pr}$ in $\gm_{\sr \Iscr\pr}=\pi(\gm_{\sr \Iscr}).$
The implication (i) $\Longrightarrow$ (ii) in 2.19 follows from the
\proclaim Lemma. With the above hypotheses
$$I(\ov\Vscr_T)=S(\gm_{\sr \Iscr}^-)I(\ov\Vscr_{T\pr}).$$
\par
\Pf
From the commuting diagram
$$\matrix{S(\pi(\gm_{\sr \Iscr}^-))&{\buildrel \pi^*\over\hookrightarrow }&S(\gm_{\sr \Iscr}^-)\cr
                  \downarrow &  &      \downarrow          \cr
            R[\ov\Vscr_{T\pr}]&{\buildrel \pi^*\over\hookrightarrow }&R[\ov\Vscr_T]\cr}$$
we obtain the inclusion $\supset.$ On the other hand 
$$S(\gm_{\sr \Iscr}^-)\slash IS(\gm^-_{\sr \Iscr})\cong S(\pi(\gm_{\sr \Iscr}^-))\slash I\otimes S(\ker \pi),$$
so $IS(\gm_{\sr \Iscr}^-)$ is prime. Equality of codimensions finishes the proof.
\QED
\parno
{\bf 3.14}\ \ \ Let $V$ be a vector space of dimension $n<\infty.$
Let $\{v_i\}_{i=1}^n$ be a basis of $V$ and set $v_{\sr 0}=v_{n+\sr 1}=0.$
Define $e,f\in \End V,$ by $ev_i=v_{i-\sr 1},\ fv_i=v_{i+\sr 1},\ \forall i.$
Let $r,s$ be integers $>0$ and set $A_{r,s}^n(t)=\Det (te^s+f^r).$
\proclaim Lemma. For all $r,s,n>0$ one has
$$A_{r,s}^n(t)=\cases{((-1)^{r+s}t^r)^h & if $h(r+s)=n,$\cr
                      0 & if $r+s$ does not divide $n.$\cr}$$
\par
\Pf
If $r+s>n,$ then there are only zeros in the $n-s+1$-th row, so we can 
suppose $r+s\leq n.$ Then the first $r$ rows have exactly one entry, namely $t.$ Similarly 
the first $s$ columns have exactly one entry, namely 1. Developing these rows and columns gives
$$A_{r,s}^n(t)=(-1)^{r+s}t^rA_{r,s}^{n-r-s}(t),$$
and hence the required result.
\QED
\parno
{\bf 3.15}\ \ \ The irreducibility of $m_{l_{\sr \Iscr}}$ is established by induction on $n.$
In view of 3.13 it is therefore enough to consider
only the situation described in 2.19 (i). Precisely we show the
\proclaim Proposition. Take $\Iscr \in\Pi$ and $T_\Iscr=(C_1,C_2,\ldots,C_l)$ the chain 
form of $T_\Iscr.$ Set $c_i=|C_i|.$ Assume $c_{\sr 1}=c_l=:c$ and 
$c_j\ne c,\ \forall j\ :\ 1<j<l.$ Then $m_{l_{\sr \Iscr}}$ is irreducible.
\par
\Pf
Otherwise by 2.17, 3.12 and 3.13, there is a subset $\Iscr\pr$ of $\Pi$
obtained from $\Iscr$ exactly as in 2.19 (ii) such that the corresponding $f\pr:=m_{l_{\Iscr\pr}}$
is an irreducible factor of 
\break
$f=m_{l_{\sr \Iscr}}.$ We obtain a contradiction by showing there exists a 
point $x\in\gm$ such that 
$f(x)=(-1)^c,\ f\pr(x)=0.$
\par
Set $a=|T_\Iscr^{c+1,k}|$ and 
$$<T_\Iscr^{1,c}>=\{p_1,p_2,\ldots,p_{n-a}\},\quad <T_\Iscr^{c+1,k}>=\{q_1,q_2,\ldots,q_a\},$$
written in increasing order. Define $x\in \gm$ by
$$x_{r,s}=\cases{1,& if $r=p_i,$ for $i:\ c<i\leq n-a$ and $s=p_{i-c}$\cr
                 0,& otherwise. \cr}$$
Set $M:=M_\Iscr^c(t).$ Let $b$ be an integer $0<b\leq n-c$ and let $\Delta_b$ denote
the set of subsets of $\{c+1,\ldots,n\}$ of cardinality $b$. Given 
$\sigma=\{r_{\sr 1},r_{\sr 2},\ldots,r_b\}\in\Delta_b,$ let $M^\sigma$ be obtained from $M$ by deleting the
$(r_i-c)-$th rows and $r_i-$th columns where $i=1,2,\ldots,b.$ Let $M_o^\sigma$ be the 
evaluation of $M^\sigma$ at $t=0.$
Developing $M$ gives the term $(-1)^c t^b M^\sigma$ and moreover
$$m_{n-c-b}=(-1)^c\sum\limits_{\sigma\in\Delta_b}M_o^\sigma.\eqno{(*)}$$
\par
Recall that $l_{\sr \Iscr}=n-c-a.$ For $\sigma\in\Delta_a$ one has
$$M_o^\sigma(x)=\cases{1, & if $\sigma=\{q_1,q_2,\ldots,q_a\},$\cr
                       0,& otherwise.\cr}$$
\par
Consequently
\item{(i)} $f(x)=(-1)^c.$
\par
\bigskip
In the notation of theorem 2.19 (ii) set $c\pr:=c_i=c_j$ and $M\pr=M^{c\pr}_{\Iscr\pr}(t).$ Then
\item{(ii)} $f\pr(x)=0,$ if $c\pr>c.$
\par
Indeed in this case the last row of $M\pr$ is indexed by the largest integer in $C_j.$
By the hypothesis this is an entry of $T_\Iscr^{c+1,k}.$ Hence even $M\pr(x)=0.$
\par
\bigskip
Define $u,v$ as in theorem 2.19 (ii) and set $n\pr=v-u+1.$
Set 
\break
$\sigma\pr=<T_\Iscr^{c+1,k}>\cap\{u,u+1,\ldots,v\}$ and $d\pr=|\sigma\pr|.$ In the notation of
3.14 we have 
\item{(iii)} $M\pr(x)=(-1)^{c\pr}t^{d\pr}A^{n\pr}_{c-c\pr,c\pr}(t),$ if $c\pr<c.$
\par
This is obtained through a development of $M\pr$ similar to that for $M$ and noting that 
$M_o^{\prime\sigma}(x)=0$ unless $\sigma\sp\sigma\pr,$ so that 
$M\pr(x)=(-1)^{c\pr}t^{d\pr}M^{\prime\sigma\pr}(x)=(-1)^{c\pr}t^{d\pr}A^{n\pr}_{c-c\pr,c\pr}(t).$
\par
\bigskip
Let us show finally that
\item{(iv)} $f\pr(x)=0,$ if $c\pr<c.$ 
\par
By 3.14 it suffices to show that $h(c-c\pr)+d\pr>d_{\Iscr\pr}$ when $n\pr=hc.$
\par
Let $C_{m_1},C_{m_2},\ldots,C_{m_r}$ be the chains of length $>c\pr$ between $C_i$ and $C_j.$
One has $d_{\Iscr\pr}=\sum\limits_{s=1}^r(c_{m_s}-c\pr)\leq d\pr+r(c-c\pr)$ and 
so the assertion holds if $r<h.$ On the other hand $hc=n\pr\geq (r+2)c\pr+d_{\Iscr\pr}.$
Thus 
$$h(c-c\pr)+d\pr\geq (r+2)c\pr+d_{\Iscr\pr}>d_{\Iscr\pr}\quad {\rm if}\quad r\geq h.$$
This proves (iv).
\par
Finally the hypothesis itself excludes $c=c\pr$ and so the proposition follows from (i), 
(ii) and (iv).
\QED
\parno
{\bf Remark.} The proof of theorem 2.19 is now complete.
\parno
{\bf 3.16}\ \ \ Adopt the hypotheses of 2.19 (i). Let $r$ be the number of chains $C_i$
of length $>c$ and $s+2c$ the sum of the lengths of the remaining chains.
\proclaim Proposition. The generator $f_{\sr \Iscr}(l)$ of $I(\ov\Vscr_{T_\Iscr(l)})$ in $S(\gm^-)$
is the highest common divisor of the (non-vanishing) $(s+c)\times(s+c)$ minors of $\M^{r+1}_\Iscr.$
\par
\Pf
Under the hypotheses of 2.19 (i), it follows from 2.18 (in notation of 2.18) that 
$\gm_{\sr \Iscr}\cap\Oscr_\Iscr(c)=\Vscr_{T_\Iscr(l)}.$ By 2.4 (i) the nilpotent orbit to which
an element $x$ belongs is determined by $\{\rk x^i\ :\ i=1,2,\ldots,n\}.$
\par
Set $\l=\sh T_\Iscr.$ Since $n-\sum_{j>r+1}\l^*_j=s+c,$ it follows that
$\rk x^{r+1}=s+c,$ for all $x\in\gm\cap\Oscr_\Iscr.$
\par
Now $\sh T_\Iscr(l)$ is obtained from $\l$ by lowering a box from $c-$th to the 
$(c+1)-$th row. Under the hypotheses of 2.19 (i) this takes it from the $(r+2)-$th
to the $(r+1)-$th column.
\par
From the three paragraphs above we conclude that
$$\gm_{\sr \Iscr}\cap\Oscr_\Iscr(c)=\{x\in\gm\ |\ \rk x^{r+1}=s+c-1\}$$
Thus this hypersurface in $\gm_{\sr \Iscr}$ is the set of common zeros of the $(s+c)\times (s+c)$
minors of $\M_\Iscr^{r+1}.$ Hence $f_{\sr \Iscr}(l)$ is exactly their largest common divisor.
\QED
\parno
{\bf 3.17}\ \ \  It is of course not too easy to explicitly determine $f_{\sr \Iscr}(l)$ through 3.16.
Moreover these power rank conditions (i.e. 2.4 (i)) are in general insufficient to obtain ideals 
of definition of orbital variety closures. This was already 
observed by van Leeuwen [vanL, \S 8] and also 
resulted independently from [M5, \S 4.3] because otherwise the chain order defined in [M5, 4.3.1]
would coincide with the geometric order which fails in $\gs\gl_7$  as it is shown in [M5, 4.3.6].
\par
Combining 3.16 and 2.19 (i) (which implies that $f_{\sr \Iscr}(l)=m_{l_{\sr \Iscr}}$) gives a remarkable combinatorial
fact, namely 3.16, about generic matrices. Here we note that $\deg m_{l_{\sr \Iscr}}=(r+1)c+s$ whilst the minors in the 
conclusion of 3.16 have degree $(r+1)(s+c).$ Thus these degrees coincide if and only if $rs=0.$
Of course when degrees coincide $f_{\sr \Iscr}(l)$ is proportional to any non-zero $(s+c)\times (s+c)$
minor of $\M_\Iscr^{r+1}.$ Using this Benlolo and Sanderson were able to prove their conjecture, 
that is Theorem 2.19 (i), under hypothesis $r=0$ or the hypothesis $c=1$ (which forces $s=0$).
These cases are relatively easy up to proving irreducibility
for which they developed a special trick. 
A purely combinatorial proof of the case $s=0$ already appears to be rather difficult.
\bigskip
\parno
{\bf 4. Strong quantization of hypersurface orbital varieties} 
\parno
{\bf 4.1}\ \ \ We begin with a combinatorial lemma. Let $\Re$ be the real field.
\proclaim Lemma. Let  $c,c_{\sr 2},\ldots,c_l$ be positive integers with
$c=c_l.$
Suppose that the system of inequalities for $\{b_i\}_{i=2}^l$ defined by
\item{1)} $1+\sum\limits_{i=2}^s b_i\leq \sum\limits_{i=2}^{s-1} c_i+
\max(c,c_s)\ :\ s=2,\ldots,l-1.$
\item{2)} $1+\sum\limits_{i=s+1}^l b_i\leq 
\sum\limits_{i=s+1}^{l-1} c_i+\max(c,c_s)\ :\ s=2,3,\ldots,l-1.$
\item{3)} $1+\sum\limits_{i=2}^{l} b_i=\sum\limits_{i=2}^{l}c_i$\hfill
\parno 
has a solution $b_i\in\Re,\ \forall i.$ Then 
$c_s\ne c,\ \forall s:\ 2\leq s<l.$ Moreover, if the latter holds one can find a solution with
$b_i\in {\Na},\ \forall i.$
\par
\Pf
Assume that the system has a solution $b_i\in\Re,\ \forall i.$
Adding 1) and 2) gives
$$2+\sum\limits_{i=2}^l b_i\leq \sum\limits_{i=2}^l c_i +2\max(c,c_s)-c-c_s$$
Then in view of 3) we obtain $1\leq 2\max(c,c_s)-c-c_s,$ that is 
$c\ne c_s$ for any $s:\ s<k.$
\par
Conversely equality in 1) and 3) gives
$$b_i=\cases{\max(c,c_{\sr 2})-1, & if $i=2$\cr
             c_{i-\sr 1}+\max(c,c_i)-\max(c,c_{i-\sr 1}), & if $2<i\leq l$\cr}\eqno{(*)}$$
These imply that
$$\sum\limits_{i=s+1}^l b_i=\sum\limits_{i=s}^l c_i-\max(c,c_s)=
\sum\limits_{i=s+1}^{l-1} c_i+\max(c,c_s)+((c+c_s)-2\max(c,c_s))$$
Then 2) holds if and only if $c+c_s\leq 2\max(c,c_s)-1$ for all $s=2,\ldots,l-1.$ Obviously
this is equivalent to $c_s\ne c$ for any $s=2,\ldots,l-1.$
Note that by $(*)$ these inequalities provide a solution with $b_i\in \Na$ for all $i.$
\QED
\parno
{\bf 4.2}\ \ \ Recall the notation of 2.7.
Given $\Iscr\st \Pi$ let $T_\Iscr=(C_1,\ldots,C_l)$ be the corresponding 
tableau. Assume that $c_{\sr 1}=c_l=c.$
One has $\Pi\pr:=\Pi\setminus \Iscr=\{\a_{\s_i}\ :\ i=1,2,\ldots,l-1\}.$
For $i=1,2,\ldots,c$ we set 
$$\b_i=\sum\limits_{j=i}^{n-1-i}\a_j.$$
Set $\b=\b_1.$ Recall that the BS element $m_{l_{\sr \Iscr}}$ defined in 2.18 has weight 
$\wt(m_{l_{\sr \Iscr}})=-\sum\limits_{i=1}^c \b_i.$
\par
Call $\nu\in\gh^*,\ \Iscr-$regular if $(\a ,\nu+\rho)\ne 0,\ \forall \a\in R_\Iscr.$
Set $w.\nu:=w(\nu+\rho)-\rho.$
Set 
$$\eqalign{\Sscr_\nu&=\{(m,\gamma)\in \Na^+\times R^+\setminus R_\Iscr^+\ |\ (\g^\vee ,\nu+\rho)=m\}\cr
           \Sscr_\nu^o&=\{(m,\gamma)\in\Sscr_\nu\ |\ s_\g.\nu+\rho\ {\rm is}\ \Iscr\ {\rm regular}\}\cr}$$
Take  $c_i\ne c \ i=2,3,\ldots,l-1$ in 4.1 and recall our normalization, making 
$(\a,\a)=2,\ \forall \a\in R.$
Define $\mu\in P(\Pi),$ by $(\mu+\rho,\a)=1,$ for $\a\in\Iscr$ and $(\mu+\rho,\a_{s_i})=-(b_{i+1}-1)$
for $i=1,2,\ldots,l-1$ and where $\{b_i\}_{i=2}^l$ are given by the conclusion of 4.1.
\proclaim Proposition. $\Sscr_\mu^o=\{(c,\b)\}.$
\par
\Pf
One has 
$$\sum\limits_{\a\in\Iscr}(\mu+\rho,\a)=\sum\limits_{i=1}^l(c_i-1)=\sum\limits_{i=2}^l(c_i-1)+c-1=
\sum\limits_{i=2}^l(b_i-1)+c=-\sum_{\a\in \Pi\pr}(\mu+\rho,\a)+c$$
by 4.1 (3) and so 
$$(\mu+\rho,\b^\vee)=(\mu+\rho,\b)=\sum\limits_{\a\in\Iscr}(\mu+\rho,\a)+\sum_{\a\in \Pi\pr}(\mu+\rho,\a)=c.$$ 
Thus $(c,\b)\in\Sscr_\mu.$
\par
Observe that if $(c\pr,\gamma)\in \Sscr_\mu\ :\ \gamma\in R^+\setminus R^+_\Iscr,\ c\pr\in \Na,$
then $(c\pr,\gamma)\in\Sscr_\mu^o$ means that 
$$s_\gamma(\mu+\rho)=\mu+\rho -(\gamma^\vee ,\mu+\rho)\gamma=\mu+\rho-c\pr\gamma$$
does not vanish on any $\d\in R^+_\Iscr.$
Note that for $\d\in R^+_\Iscr$ and $\g\in R^+\setminus R_\Iscr^+$ one has
$$(\d,\g)=\cases{0 &if $\g-\d\not\in R^+$\cr
                     1 & if $\g-\d\in R^+$\cr}$$
Thus vanishing on some $\d\in R^+_\Iscr$ means that $\g-\d$ is a root and $(\d,\mu+\rho)=c\pr$
or simply that $\g-\d$ is a root such that $(\mu+\rho,\g-\d)=(\mu+\rho,\g)-(\mu+\rho,\d)=0.$
\par
Note that $\b-\d$ is a root for $\d\in R^+_\Iscr$ iff
$\d=\a_1+\cdots+\a_i$ or $\d=\a_{n-1}+\cdots+\a_{n-i}$  where $i<c.$ Thus $(\mu+\rho,\b-\d)>0$
and so $(c,\b)\in\Sscr_\mu^o.$
\par 
Consider $\g=\a_i+\a_{i+1}+\cdots+\a_j$ for $i<c$ and $j>n-c.$ One has $(\mu+\rho,\g)=c-(i-1)-(n-j-1)=:c\pr.$
Then unless $\g=\b$ 
one can find $\d\in R^+_\Iscr$ starting at $\a_i$ or ending in $\a_j$ 
such that $(\mu+\rho,\g-\d)=0.$ We conclude that $(c\pr,\g)\not\in\Sscr_\mu^o.$
\par
Set $\d_i^-=\a_1+\a_2+\cdots+\a_{\s_{i-1}}$ for $2\leq i\leq l$ and 
$\d_j^+=\a_{\s_j}+\a_{\s_j+1}+\cdots+\a_{n-1}$
for $1\leq j\leq l-1.$ One has 
$$\eqalign{
(\mu+\rho,\d_s^-)&=\sum\limits_{i=1}^{s-1}(c_i-1)+\sum_{i=2}^s -(b_i-1)=c+
                    \sum\limits_{i=2}^{s-1} c_i-\sum_{i=2}^s b_i\geq c+1-\max(c,c_s)\cr
(\mu+\rho,\d^+_s)&=\sum\limits_{i=s+1}^l(c_i-1)+\sum\limits_{i=s+1}^l-(b_i-1)\cr
                 &=
\sum\limits_{i=s+1}^{l-1}c_i-\sum\limits_{i=s+1}^l b_i+c\geq c-\max(c,c_s)+1\cr}
$$
Now consider $\g_{s,l}=\b-\d_s^-$ with $2\leq s\leq l-1$
Then $(\mu+\rho,\g_{s,l})\leq \max(c,c_s)-1.$ Thus it follows that $(c\pr,\g_{s,l})\not\in \Sscr_\mu^o.$
Indeed if $c_s> c$ we can subtract $\d=\a_{\s_{s-1}+1}+\cdots+\a_i$ for some $i\ :\ \s_{s-1}+1\leq i<\s_s$  
from $\g_{s,l}$ so that $(\mu+\rho,\g_{s,l}-\d)=0$
and if $c_s< c$ we can subtract $\d=\a_{n-1}+\cdots \a_i$ for some $i\ :\ \s_{l-1}<i\leq n-1$ 
from $\g_{s,l}$ so that $(\mu+\rho,\g_{s,l}-\d)=0.$
A similar assertion holds if replace $\g_{s,l}$ by $\g=\g_{s,l}-(\a_{\s_{s-1}+1}+\cdots+\a_i)-(\a_{n-1}+\cdots \a_j)$
for any $i\ :\ \s_{s-1}+1\leq i<\s_s$ and $j\ :\ \s_{l-1}<j\leq n-1,$ that is to say when some
of these roots have already been subtracted.
\par
A similar conclusion holds for $\g_{1,s}=\b-\d_s^+$ with $2\leq s\leq l-1$
and for $\g= \g_{1,s}-(\a_1+\cdots+\a_i)-(\a_{\s_s-1}+\cdots \a_j)$ for any
$i\ :\ 1\leq i<c$ and $j\ :\ \s_{s-1}< j\leq \s_s-1$ (by the obvious symmetry).
\par
Finally consider $\g_{i,j}=\b-\d_i^- -\d_j^+,$ with $2\leq i<j\leq l-1.$
We have 
$$(\mu+\rho,\g_{i,j})\leq \max (c,c_i)+\max(c,c_j)-c-2.$$
Then a similar conclusion holds in this case also with slightly stronger reason.
\QED
\parno
{\bf Remark.} One may also check that the inequalities 1) and 2) of 4.1 are necessary for the conclusion of the 
proposition to hold with $\mu+\rho$ satisfying $(\mu+\rho,\a)=1,\ \forall \a\in \Iscr$ and $(\mu+\rho,\b)=c.$
Consequently the proposition fails when $c=c_i,$ for some $i:\ 1<i<l.$
\parno
{\bf 4.3}\ \ \ Take  $\mu$ as in 4.2 and recall the definition of $\b_i$ from 2.20.
Recall that $\g_c=\sum\limits_{i=1}^c \b_i.$ Then $(\g_c,\a)=0,\ \forall \a\in \Iscr.$
Moreover $\mu+\rho-\g_c$ is the unique $\Iscr$ dominant element of $W_\Iscr(\mu+\rho-c\b)=
W_\Iscr(s_\b(\mu+\rho)).$ 
\par
Let $\{M_\Iscr^i(\mu)\}_{i=0}^\infty$ be the Jantzen filtration of $M_\Iscr(\mu).$
Given $(m,\g)\in\Sscr^o_\mu,$ let $\omega(s_\g.\mu)$ be the unique $\Pi\pr$ dominant element in
$W_\Iscr.(s_\g.\mu)$ and $\sign(s_\g.\lambda):=(-1)^{l(w)}$ where $w\in W_\Iscr$ is the unique element
satisfying $w.s_\g.\mu=\omega(s_\g.\mu).$
The Jantzen sum formula combined with 4.2 gives
$$\sum\limits_{i=1}^\infty \ch M_\Iscr^i(\mu)=\sum\limits_{(m,\g)\in \Sscr^o_\mu}\sign(s_\g.\mu)\ch M_\Iscr(\omega(s_\g.\mu))
=\ch M_\Iscr(\mu-\g_c), \eqno{(*)}$$
by the above noting that $\sign(s_\b.\mu)$ is necessarily positive because there
is only this term on the right hand side. One may recall that the right hand side of the Jantzen sum 
formula obtains from the zero of the Shapovalov determinants. Then the appearance of just one term
on the right hand side means that the Shapovalov determinant is non-zero on $M_\Iscr(\mu)_\nu\ :\ \nu>\omega(s_\b.\mu)$
and vanishes on $M_\Iscr(\mu)_{\omega(s_\b.\mu)}.$ 
\par
(Equivalently, but more directly by Jantzen [Ja, Satz 2] the Shapovalov determinant on the
$\nu$ weight subspace of $M_{\Iscr(\mu)}$ takes form
$$\prod\limits_{m=1}^\infty\prod\limits_{\g\in R^+\setminus R^+_\Iscr}((\g,\mu+\rho)-m)^{\chi\pr_\nu(\mu-m\g)}$$
where 
$$\chi\pr_\nu(\lambda)=\sum\limits_{w\in W_\Iscr}(-1)^{l(w)}\dim M(w.\lambda)_\nu \eqno{(*)}$$
This has a zero if and only if some $(m,\g)\in \Sscr_\mu$ and the sum of the corresponding exponents is
strictly positive. Now by $(*)$, $\chi\pr_\nu(\lambda)\ne 0$ for some $\nu$ if and only if $\lambda+\rho$
is $\Iscr$ regular, so we further require $(\mu,\g)\in\Sscr^o_\mu.$
Moreover in this case $\chi\pr_\nu(\lambda)$ equals $\sign(\lambda)\dim M_\Iscr(\omega(\lambda))_\nu.$
Since in addition $\Sscr^o_\mu=\{(m,\b)\}$ in the present case we must have 
$\chi\pr_\nu(s_\b.\mu)\geq 0$ for all $\nu$ and of course no cancellations occur.
Moreover from $(*)$ we obtain
$$\chi\pr_\nu(s_\b.\mu)=\cases{0, & if $\nu>\omega(s_\b.\mu),$\cr
                               1, & if $\nu=\omega(s_\b.\mu).$\cr}$$
which gives the required assertion.)
\par
This forces there to be a highest weight vector
in $M_\Iscr(\mu)$ of weight $\omega(s_\b.\mu)=\mu-\g_c.$ Since $(\a,\mu)=0$
for all $\a\in\Iscr$ it follows that $M_\Iscr(\mu)$ is a rank 1 free $U(\gm_{\sr \Iscr}^-)$
module. Again $(\a,\mu-\g_c)=0,$ for all $\a\in\Iscr.$ Consequently this highest weight vector
must generate a submodule of $M_\Iscr(\mu)$ isomorphic to $M_\Iscr(\mu-\g_c)$ necessarily contained in the
maximal submodule $M_\Iscr^1(\mu)$ of $M_\Iscr(\mu).$ By $(*)$ we then conclude that 
$M_\Iscr^1(\mu)=M_\Iscr(\mu-\g_c)$ and so $V(\mu)=M_\Iscr(\mu)/M_\Iscr(\mu-\g_c)$ is simple.   
\par
By [M4] the associated variety of $V(\mu)$ is irreducible and hence an orbital variety $\Vscr.$
\par
Set $p_{\sr \Iscr}=\prod\limits_{\a\in R^+_\Iscr}\a.$
Since
$$\ch V(\mu)=(\ch M_\Iscr(\mu))(1-e^{-\g_c})=\ch S(\gm^-_{\sr \Iscr})e^\mu(1-e^{-\g_c}),$$
it follows that the characteristic polynomial of $\Vscr$ is just $\g_c p_{\sr \Iscr}.$
\par
Let $f_{\sr \Iscr}(l)$ be the BS element constructed in 2.18 and $\Vscr_{T_\Iscr(l)}$ the corresponding 
hypersurface orbital variety. By 3.12 the zero variety $\Vscr(f_{\sr \Iscr}(l))$ 
of $f_{\sr \Iscr}(l)$ in $\gm_{\sr \Iscr}^-$ is a union of orbital variety closures of codimension 1, one of which 
is $\Vscr_{T_\Iscr(l)}.$ Yet $\wt f_{\sr \Iscr}(l)=-\g_c$ and so the characteristic polynomial of $\Vscr(f_{\sr \Iscr}(l))$ 
is also $\g_c p_{\sr \Iscr}.$ By the linear independence of characteristic polynomials of orbital varieties [J1] we
conclude that $\Vscr=\Vscr(f_{\sr \Iscr}(l))=\Vscr_{T_\Iscr(l)}.$ This proves that $f_{\sr \Iscr}(l)$ is
irreducible. Finally $R[\ov\Vscr_{T_\Iscr(l)}]=S(\gm_{\sr \Iscr}^-)/S(\gm_{\sr \Iscr}^-)f_{\sr \Iscr}(l),$ so
$\ch R[\ov\Vscr_{T_\Iscr(l)}]=\ch S(\gm_{\sr \Iscr}^-)(1-e^{-\g_c}).$ Comparison with $\ch V(\mu)$ which is above 
expression up to a shift, shows
that $V(\mu)$ is a strong quantization of $\Vscr_{T_\Iscr(l)}.$
We have proved the
\proclaim Theorem. Every hypersurface orbital variety in $\gs\gl_n$ admits a
strong quantization being a simple highest weight module with integral
highest weight given by the above procedure.
\par
\parno
{\bf 4.4}\ \ \ Let us consider the case where we take equality in 1) of 4.1. Using 3)
we obtain 
$$\sum\limits_{i=2}^s b_i=\sum\limits_{i=2}^{s-1}c_i+\max(c,c_s)-1,\quad 
\sum\limits_{i=s+1}^l b_i=\sum\limits_{i=s}^l c_i-\max(c,c_s)$$
Then in the notation of 
4.2 for $2\leq s\leq l-1$ we obtain equalities
$$\eqalign{
(\mu+\rho,\d_s^-)&=\sum\limits_{i=1}^{s-1}(c_i-1)+\sum_{i=2}^s -(b_i-1)=c-\max(c,c_s)+1 \cr
(\mu+\rho,\d^+_s)&=\sum\limits_{i=s+1}^l(c_i-1)+\sum\limits_{i=s+1}^l-(b_i-1)= \max(c,c_s)-c_s\cr}$$
Consequently $(\mu+\rho,\g_{s,l})= \max(c,c_s)-1$ for $2\leq s\leq l-1$ and 
$(\mu+\rho,\g_{1,s})=c+c_s-\max(c,c_s)$ for $2\leq s< l.$
Finally for any $\g_{s,t}$ with $1< s<t< l$ one has
$$\eqalignno{(\mu+\rho,\g_{s,t})= &(\mu+\rho,\g_{1,t})-(\mu+\rho,\d_s^-)
=c+c_t-\max(c,c_t)-(c-\max(c,c_s)+1)\cr
=&(c_t-1)-\max(c,c_t)+\max(c,c_s).&(*)\cr}$$
\par
Now let $\g_{i,j}^-$ denote $\g_{i,j}$ with the $(c_i-1)$ roots of the $i$-th column and $(c_j-1)$ of the $j$-th column removed
that is 
$$\g_{i,j}^-=\a_{\s_{i+1}}+\cdots+\a_{\s_{j-1}}.$$
Then
$$(\mu+\rho,\g_{1,t}^-)=c+c_t-\max(c,c_t)-(c-1)-(c_t-1)=2-\max(c,c_t).$$
Yet to $\g_{1,t}^-$ we can always add either the roots from the first column
or from $t$-th column to obtain some $\g:=\d+\g^-_{1,t}\ :\ \d\in R^+_\Iscr.$
Such $\g$ for which $(\mu+\rho,\g)$ is maximum satisfies
$$(\mu+\rho,\g)=2-\max(c,c_t)+\max(c,c_t)-1=1>0.\eqno{(**)}$$
Similarly for $2\leq s<t\leq l$ by $(*)$ one has
$$(\mu+\rho,\g_{s,t}^-)=\max(c,c_s)-\max(c,c_t)-(c_s-1).$$
To $\g_{s,t}^-$ we can add either the roots from the $s$-th column or from
$t$-th column. Again in the best case for such $\g$ one has
$$(\mu+\rho,\g)=\max(c_s,c_t)+\max(c,c_s)-\max(c,c_t)-c_s\geq 0.\eqno{(***)}$$
Let us define
$$\check\Sscr_\mu=\{(m,\b)\in{-\Na^+}\times R^+\setminus R_{\Iscr}^+\ |\ (\b^\vee,\mu+\rho)=m\}$$
and
$$\check\Sscr^o_\mu=\{(m,\b)\in \check\Sscr_\mu\ |\ s_\b(\mu+\rho) {\rm\ is\ } \Iscr {\rm\ regular\ }\}.$$
\proclaim Lemma. With the above choice of $\mu$ one has $\check\Sscr^o_\mu=\emptyset.$
\par
\Pf
As before given $\g\in \check\Sscr_\mu$ we must find $\d\in R^+_\Iscr$ such that $0=(\d,s_\g(\mu+\rho))=
(\d,\mu+\rho+(\mu+\rho,\g)\g).$ It is enough to show that $\g+\d$ is a root vanishing on
$\mu+\rho.$ For this notice for $\g\in \check\Sscr_\mu$ we must first have 
$r:=(\g,\mu+\rho)\in -{\Na}^+.$ If for example 
$\g=\g_{s,t}^-\ :\ 1\leq s<t\leq l,$ then $(**),\ (***)$ show that there exists $\d\in R^+_\Iscr$
such that $\g+\d$ is a root and $(\g+\d,\mu+\rho)$ takes all possible integer values from $r$
to an integer $\geq 0.$ The general case is similar (and easier). 
\QED
\parno
{\bf 4.5}\ \ \ Let $\bO_\Iscr$ be the subcategory of $\bO$ in which the Levi factor defined by $\Iscr$
acts finitely. By [JLT, 9.6] 
\proclaim Proposition. Given $\check\Sscr^o_\mu=\emptyset,$ then $M_\Iscr(\mu)$ is projective in $\bO_\Iscr.$
\par
\proclaim {4.6\ \ \ Theorem}. Define  $\mu$ as in 4.4. Then $\Ann V(\mu)\in \Max U(\gg).$
\par
\Pf 
Since $\gg$ is of type $A$ the natural map $U(\gg)\ra F(M_\Iscr(\mu),M_\Iscr(\mu))$ is surjective.
Set $F^\mu=F(M_\Iscr(\mu),M_\Iscr(\mu))=U(\gg)\slash \Ann M_\Iscr(\mu).$ Let $P$ be a maximal ideal of $F^\mu.$
By [JLT, 10.9], 4.4 and 4.5 we have
$$P=\Ann_{F^\mu}(M_\Iscr(\mu)\slash PM_\Iscr(\mu)).$$
Hence $PM_\Iscr(\mu)\subsetneq M_\Iscr(\mu).$ Yet $M_\Iscr(\mu-\g_c)$ is the maximal submodule of
$M_\Iscr(\mu)$ as a $U(\gg)$ module. This forces $P\subset \Ann V(\mu)$ and hence equality.
\QED
\parno
{\bf 4.7}\ \ \ Take $c=c_{\sr 1},c_{\sr 2},\ldots,c_l=c.$ Proceeding as in 4.2 one may rather easily 
show that there exists $\mu\in\gh^*$ satisfying $(\mu+\rho,\a)=1,$ for all $\a\in\Iscr$ and such that
$\Sscr^o_\mu=\{(c,\b)\}.$
Then as in 4.3, one checks that $V(\mu)$ is a strong quantization of $\Vscr(f_{\sr \Iscr}(l)).$
Here we do not need $c_i\ne c,$ for $i\ :\ 1<i<l,$ however the resulting $\mu$ will not be integrable, nor
will the associated variety of $V(\mu)$ be irreducible. However it will be contained in the
associated variety of $\Ann V(\mu)$ which by Borho-Kraft [BK] is just the closure of a nilpotent orbit
specifically $\Oscr_{T_\Iscr(l)}:=\bG\Vscr_{T_\Iscr(l)}$ in this case. 
Consequently the associated variety of $V(\mu)$ is contained in $\gm\cap\Oscr_{T_\Iscr(l)}$ and we recall that the latter
is irreducible if and only if $c_i\ne c,$ for all $i\ :\ 1<i<l.$ This concludes the third proof of the
irreducibility of the BS elements. 
\par
\bigskip
\ctr{{\bf Appendix: Index of Notation}}
Symbols appearing frequently are given below in order of appearance.
\item {1.1} $\gg$
\item {1.2} $\gn^-,\ \gn^+,\  \gh,\ \gb,\ \bG,\ \bBo,\ S(\cdot),\ U(\cdot),\ \Oscr,\ \Vscr,\ I(\ov\Vscr)$
\item {1.3} $\l,\ V(\l),\ \Fscr$
\item {1.5} $\gm^+,\ \gm^-,\  \gp,\ \Mscr(t),\ M(t)$
\item {1.7} $\bP$
\item {2.1} $\gn,\ \bV,\ W,\ R,\ R^+,\ \Pi,\ X_\a,\ \gn\cap^w\gn,\ \Vscr_w$
\item {2.2} $\gg,\ \bG,\ \gn,\ \gn^-,\ \bBo,\ \ e_{i,j},\ \a_{i,j},\ \a_i,\ s_\a,\ s_i$
\item {2.3} $\bS_n,\ [a_{\sr 1},\ldots,a_n],\ p_w(i),\ S(w)$
\item {2.4} $P(n),\ \l,\ \l^*,\ k,\ l,\ D_\l$
\item {2.6} $\bT_n,\ [a\pr, a\prpr],\ \sh T,\ \sigma(T),\ \widehat \cdot,\ <\cdot>,\ |\cdot|,\ Q(w),\ \Vscr_w,\ T_\Vscr,\ \Vscr_1>\Vscr_2,\ T_1>T_2$
\item {2.7} $T_j^i,\ r_{\sr T}(\cdot),\ c_{\sr T}(\cdot),\ T^i,\ T_j,\ h(T^i_j),\ T^{i,j},\ (T,S),\ {T\choose S}$
\item {2.8} $(R+j),\ w_r(T),\ w_c(T),\ (T\da b)$
\item {2.9} $(R\ua j),\ (R-a),\ (C\la j),\ (T-T^i_j)$
\item {2.10} $\tau(w),\ \tau(T),\ \tau(\bP),\ \tau(\gp),\ \tau(\Vscr),\ \bP_T,\ \bP_\Vscr,\ \gp_{\sr T},\ \gp_{\sr \Vscr}$ 
\item {2.11} $\Iscr,\ \bP_\Iscr,\ \bM_\Iscr,\ \bL_\Iscr,\ \gp_{\sr \Iscr},\ \gm_{\sr \Iscr},\ \gl_{\sr \Iscr},\ \ell(w),\ W_\Iscr,\ 
              w_{\sr \Iscr},\ \Vscr_\Iscr,\ \Oscr_\Iscr$
\item {2.12} $T_\Iscr,\ C^\Iscr_i,\ (C^\Iscr_1,C^\Iscr_2,\ldots,C_l^\Iscr)$
\item {2.14} $\bBo_\Iscr,\ \gn_{\sr \Iscr},\ W^\Iscr,\ \pi_{\sr \Iscr},\ \pi,\ \Vscr^\Iscr_{\pi(w)},\ \pi_{i,j}$
\item {2.15} $c_i,\ \s_i,\ T_\Iscr(i)$
\item {2.16} $\Iscr_n,\ \Iscr_1$
\item {2.18} $\l_{\sr \Iscr},\ \Oscr_\Iscr(c)$
\item {2.19} $f_{\sr \Iscr}(i),\ X(w),\ x_{i,j},\ \{\ ,\ \},\ \M(w),\ c,\ M_\Iscr^c(t),\ d_{\sr \Iscr},\ l_{\sr \Iscr},\ m_{l_{\sr \Iscr}}$
\item {2.20} $\b_i,\ \g_c$
\item {3.4} $f$
\item {3.5} $w-s,\ \M^{i,j}$
\item {3.8} $M$
\item {3.9} $M^{i,j}$
\item {4.2} $\b,\ \Sscr_\nu,\ \Sscr^o_\nu,\ \mu$
\item {4.4} $\check\Sscr_\mu,\ \check\Sscr^o_\mu$
\bigskip
\ctr{{\bf References}}
\item{[BK]} W. Borho and H. Kraft, {\it \"Uber Bahnen und deren Deformationen 
            bei linearen Aktionen reductiver Gruppen,} Comment. Math. Helvetici {\bf 54}
            (1979), 61-104. 
\item{[BS]} E. Benlolo and J. Sanderson, {\it On supersurface orbital varieties of} $\gs\gl
            (N,\Co)$,  J. Algebra 245 (2001), 225--246.
\item{[F]} W. Fulton, {\it Young tableaux}, LMS 35, Cambridge University Press (1997). 
\item{[H]} W. Hesselink, {\it Singularities in the nilpotent 
                       scheme of a classical group}, Trans. Am. Math. 
                       Soc. {\bf 222} (1976), 1-32.
\item{[Ja]} J.C. Jantzen, {\it Kontravariante Formen auf induzierten
            Darstellungen halbeinfacher Lie-Algebren}, Math. Ann. {\bf 226} (1977), 53-65.
\item{[J1]} A. Joseph, {\it On the variety of a highest weight module}, J. Algebra {\bf 88}
            (1984), 238-278.
\item{[J2]} A. Joseph, {\it Orbital varieties, Goldie rank polynomials and unitary highest 
            weight modules}, in B. Orsted and H. Schlichtkrull (eds.) {\it Algebraic and 
            Analytic Methods in Representation Theory}, Perspectives in Mathematics, 
            Vol.17, Academic Press, London, 1997.
\item{[J3]} A. Joseph, {\it Orbital varieties of the minimal orbit}, Ann. Ec.
            Norm. Sup. {\bf 31} (1998), 17-45.
\item{[JLT]} A. Joseph, G. Letzter and D. Todoric {\it On the 
             Kostant-Parthasarathy-Ranga Rao-Varadarajan Determinants, III. Computation 
             of the KPRV Determinants}, J. Algebra {\bf 241} (2001), 67-88.
\item{[Kn]}  D. E. Knuth, {\it The art of computer programming}, 
                       Vol.3, Addison-Wesley (1969), 49-72.
\item{[M1]} A. Melnikov, {\it On orbital variety closures in $\gs\gl_n.$ I. Induced Duflo order}, preprint
\item{[M2]} A. Melnikov, {\it On orbital variety closures in $\gs\gl_n.$ II. Descendants of a Richardson orbital variety}, 
preprint.
\item{[M3]} A. Melnikov, {\it Robinson-Schensted procedure and combinatorial
            properties of geometric order} $\gs\gl(n)$, CRAS Paris {\bf 315}, s\'erie I
            (1992), 709-714.
\item{[M4]} A. Melnikov, {\it Irreducibility of the associated variety of
            simple highest weight modules in} $\gs\gl(n)$, CRAS, Paris {\bf 316}, 
            s\'erie I (1993), 53-57.
\item{[M5]} A. Melnikov, {\it Geometric Interpretation of Robinson-Schensted
            procedure and related orders on Young tableaux}, thesis, Weizmann Institute 1992.
\item{[S]} N. Spaltenstein, {\it On the fixed point set of a unipotent element
           on the variety of Borel subgroups}, Topology {\bf 16} (1977), 203-204.
\item{[S1]} N.Spaltenstein, {\it The fixed point set of a unipotent 
            transformation on the flag manifold}, Proc. Konin. 
            Nederl. Akad. {\bf 79} (1976), 452-456.
\item{[Sch]}  M. P. Sch\"utzenberger, {\it La correspondance de Robinson}, LN in Math. {\bf 597} (1976),
                59-113.
\item{[St]} R. Steinberg, {\it An occurrence of the Robinson-Schensted
            correspondence}, J. Algebra {\bf 113} (1988), 523-528.
\item{[vanL]} M. A. A. van Leeuwen, {\it The Robinson-Schensted and Sch\"utzenberger algorithms, 
              Part II: Geometric interpretations}, CWI report AM-R9209 (1992).

\bye